\documentstyle[12pt]{article}
\def\Bbb R{{\rm \bf R}}
\def\proclaim#1{\vskip2mm{\bf #1}\em}
\def\endproclaim{\em \vskip2mm}
\def\tag#1{\eqno(#1)}
\def\gathered{\begin{array}{c}}
\def\endgathered{\end{array}}
\def\text{\mbox}
\begin{document}

\title {Reconstruction of Inclusion from Boundary Measurements}
\author{Masaru IKEHATA\footnote{
Department of Mathematics,
Faculty of Engineering,
Gunma University, Kiryu 376-8515, Japan}
}
\date{  }
\maketitle
\begin{abstract}
 We consider the problem of reconstructing  of the
boundary of an unknown inclusion together with its conductivity
from the localized Dirichlet-to-Neumann map. We give an exact
reconstruction procedure and apply the method to an inverse
boundary value problem for the system of the equations in the
theory of elasticity.

\noindent
AMS: 35R30

\end{abstract}

%\tableofcontents

\section{Introduction and statement of the results}
This paper is the sequel to \cite{IK3}.
Let $\Omega$ be a bounded domain in $\Bbb R^n$, $n=2, 3$ with connected Lipschitz
boundary.  We denote its conductivity by $\gamma$.  In what follows, unless otherwise
stated, we assume (C):
$$
\left\{
\begin{array}{l}
\displaystyle
\text{$\gamma$ is an essentially bounded real-valued function on
$\Omega$}\\
\text{and uniformly positive definite.}
\end{array}
\right.
\tag {C}
$$
We define the Dirichlet-to-Neumann map
$$\displaystyle
\Lambda_{\gamma}:H^{1/2}(\partial\Omega)\longrightarrow H^{-1/2}(\partial\Omega)
$$
by the formula
$$\displaystyle
<\Lambda_{\gamma}f, g>=\int_{\Omega}\gamma\nabla u\cdot\nabla\varphi dx,
$$
where $g\in H^{1/2}(\partial\Omega)$, $\varphi$ is any $H^1(\Omega)$ function
with $\varphi\vert_{\partial\Omega}=g$, $u$ is in $H^1(\Omega)$ and
the weak solution of the Dirichlet problem
$$\begin{array}{c}
\displaystyle
\nabla\cdot\gamma\nabla u=0\,\,\text{in}\,\Omega,\\
\\
\displaystyle
u\vert_{\partial\Omega}=f.
\end{array}
$$
Notice that $\Lambda_{\gamma}$ is a bounded linear operator.  $\Lambda_{\gamma}f$
is the electric current on $\partial\Omega$ corresponding to a voltage
potential $f$ on $\partial\Omega$.
\par
Let $\Gamma$ be a given nonempty open subset of $\partial\Omega$.
Set
$$\displaystyle
{\cal D}(\Gamma)=\{f\in H^{1/2}(\partial\Omega)\,\vert\,\text{supp}\,f\subset\Gamma\}.
$$
We call the map
$$\displaystyle
{\cal D}(\Gamma)\ni f\longmapsto \Lambda_{\gamma}f\vert_{\Gamma}
$$
the localized Dirichlet-to-Neumann map.

We assume that $\Omega$ contains an unknown inclusion ``discontinuously"
imbedded in a reference medium with known conductivity $\gamma_0$.
We denote by $D$ the shape of the inclusion and by $\gamma\vert_D$ its conductivity.

The problem is to find a reconstruction formula(procedure) of $D$ together
with $\gamma\vert_D$ by means of the localized Dirichlet-to-Neumann map.

Let us describe it more precisely.
We assume that $\gamma_0$ satisfies (C) and that
$$\begin{array}{c}
\text{$\gamma_0\in C^{0,1}(\overline\Omega)$ if $n=3$;}\\
\\
\displaystyle
\text{$\gamma_0\in C^{0,\theta}(\overline\Omega)$ with $0<\theta\le 1$ if $n=2$.}
\end{array}
$$
The pair $(D, \gamma)$ satisfies
$$\left\{
\begin{array}{l}
\displaystyle
\text{(i) $D$ is an open subset of $\Omega$
with Lipschitz boundary}\\
\text{ and satisfies $\overline D\subset\Omega$;}\\
\text{(ii) $\Omega\setminus\overline D$ is connected;}\\
\text{(iii)  $\gamma(x)=\gamma_0(x)$ for almost all $x\in \Omega\setminus\overline D$;}\\
\text{(iv) $\gamma$ satisfies the following jump condition:}\\
\forall a\in\partial D\exists\delta>0\\
\text{$\gamma-\gamma_0$ is uniformly positive definite on $D\cap B(a,\delta)$}\\
\text{or}\\
\text{$\gamma_0-\gamma$ is uniformly positive definite on $D\cap B(a,\delta)$.}
\end{array}
\right.
\tag {CI}
$$

\noindent
Our main result is

\proclaim{\noindent Main Theorem.}  Assume that $\gamma\vert_D\in W^{2,p}(D)$
with $2p>n$.  Let $\Gamma$ be a given nonempty open subset of $\partial\Omega$.
There exists a set ${\cal D}(\gamma_0,\Gamma)$ contained in ${\cal D}(\Gamma)$ such that
both $\partial D$ and $\gamma\vert_D$
can be reconstructed from the set
$$\displaystyle
\{<(\Lambda_{\gamma}-\Lambda_{\gamma_0})f,\,f>
\vert\,f\in {\cal D}(\gamma_0,\Gamma)\}.
$$
\endproclaim

{\it\noindent Remark.}
${\cal D}(\gamma_0,\Gamma)$ is independent of $\Lambda_{\gamma}$
and depends on $\Omega, \gamma_0$ and $\Gamma$.

Isakov \cite{IS} has proved the uniqueness of reconstruction
in the case when $n=3$, $\gamma_0\in C^2(\overline\Omega)$,
$\partial\Omega$ is $C^2$ and $\gamma\vert_D\in C^2(\overline D)$.
However his proof does not provide us
how to reconstruct $D$ and $\gamma\vert_D$.
Main Theorem does it.
It should be noted that Kohn-Vogelius \cite{KV} has
proved the uniqueness of reconstruction of the piecewise real analytic conductivity.
Alessandrini \cite{A} gave further uniqueness results in this direction.  Their proofs also
do not provide us any reconstruction procedure.

Nachman \cite{N2} (see also \cite{N1}) proved

\proclaim{\noindent Theorem A.}
Let $\gamma$ satisfy (C).
Assume that $\gamma\in W^{2, p}(\Omega)$ with $2p>n$.
$\gamma$ can be reconstructed from $\Lambda_{\gamma}$.
\endproclaim

Since $\gamma$ in Main Theorem has a first kind of discontinuity,
one can not immediately get Main Theorem from Theorem A.
Furthermore in Main Theorem we need only the localized
Dirichlet-to-Neumann map.   It is not clear whether one can
replace the Dirichlet-to-Neumann map in the statement of Theorem A
with the localized one.  Notice that in the case when $n=3$, Sylvester -Uhlmann \cite{SU} proved
the uniqueness of the reconstruction of $\gamma\in C^{\infty}(\overline\Omega)$
from the full knowledge of $\Lambda_{\gamma}$.

We briefly describe the steps of the proof of Main Theorem.
First in Section 2 we prove

\proclaim{\noindent Theorem B.}
Let $\Gamma$ be a given nonempty open subset of $\partial\Omega$.
Under (CI), there exists a set ${\cal D}(\gamma_0,\Gamma)$
contained in ${\cal D}(\Gamma)$ such that
$\partial D$ can be reconstructed from the set
$$\displaystyle
\{<(\Lambda_{\gamma}-\Lambda_{\gamma_0})f\,,f>\,\vert
f\in{\cal D}(\gamma_0,\Gamma)\}.
$$
\endproclaim

{\it\noindent Remark.} We do not require any regularity on
$\gamma$ inside the inclusion; if $\gamma_0\equiv1$, Theorem B is
included in \cite{IK3} and the uniqueness of the reconstruction
has been proved in \cite{IK2}; for the concrete description of
${\cal D}(\gamma_0,\Gamma)$ see Remark under Proposition 2.5 in
Section 2.

From this theorem one gets $D$.  So the next problem is to reconstruct $\gamma\vert_D$.
To do it we calculate the Dirichlet-to-Neumann map inside the inclusion.
Let us describe it more precisely.

{\bf\noindent Definition(Dirichlet-to-Neumann Map inside Inclusion).}
Consider $D$ and $\gamma$ satisfying (i), (ii) of (CI)
and (C) respectively.
We define the Dirichlet-to-Neumann map inside inclusion
$$\displaystyle
\Lambda_{\gamma^-}:H^{1/2}(\partial D)\longrightarrow H^{-1/2}(\partial D)
$$
by the formula
$$\displaystyle
<\Lambda_{\gamma^-}f, g>=\int_D\gamma\nabla u\cdot\nabla\varphi dx,
$$
where $g\in H^{1/2}(\partial D)$, $\varphi$ is any $H^1(D)$ function
with $\varphi\vert_{\partial D}=g$, $u$ is in $H^1(D)$ and
the weak solution of the Dirichlet problem
$$\begin{array}{c}
\displaystyle
\nabla\cdot\gamma\nabla u=0\,\,\text{in}\, D,\\
\\
\displaystyle
u\vert_{\partial D}=f.
\end{array}
$$

In Section 3 we prove

\proclaim{\noindent Theorem C.}  Let $\Gamma$ be a
given nonempty open subset of $\partial\Omega$. Assume that
(i)$\sim$(iii) of (CI) hold. There exists a set ${\cal D}(\gamma_0,
\Gamma, D)$ contained in ${\cal D}(\Gamma)$ such that
$\Lambda_{\gamma^{-}}$ can be calculated from the set
$$
\{<(\Lambda_{\gamma}-\Lambda_{\gamma_0})f, f>\,\vert
\, f\in{\cal D}(\gamma_0,\Gamma,D)\}.
$$
\endproclaim

{\it\noindent Remark.}  In this theorem, it is assumed that both
$\gamma_0$ and $D$ are known;  ${\cal D}(\gamma_0,\Gamma, D)$ is
independent of $\Lambda_{\gamma}$ and depends on $\Omega, \Gamma$
and $D$; notice that we do not assume any regularity on $\gamma$
inside the inclusion; for the concrete description of ${\cal
D}(\gamma_0, \Gamma, D)$ see Remark under Proposition 3.3 in
Section 3.

Thus if once we get $D$, from Theorem C  we get $\Lambda_{\gamma^{-}}$.
Then from Theorem A (Nachman's reconstruction procedure)
in the case where $\Omega=D$, $\gamma=\gamma\vert_D$
we get $\gamma\vert_D$.  For the summary of the procedure see Section 4.
\par
Theorem C has an interesting conclusion in the context of heat
conduction. Suppose you have a heat conductive body $D$ with
unknown heat conductivity $\gamma^-$.  We assume that there is no
heat sources inside $D$. You can measure the temperature
distribution $f$ on $\partial D$. But how can you measure the heat
flux $\Lambda_{\gamma^-}f$ on $\partial D$ which produces $f$ on
$\partial D$? Since $\gamma^-$ is  unknown, to do it you can not
make use of the temperature distribution at the boundary layer of
$\partial D$. A way to overcome this difficulty is to insert $D$
into a conductive body $\Omega$ with known heat conductivity
$\gamma_0$.  Then give any temperature distribution $g$ on
$\partial\Omega$ and measure the temperature at the boundary layer
of $\partial\Omega$.  Since $\gamma_0$ is known, from such data
you can calculate the heat flux on $\partial\Omega$ which produces
$g$. If you do this procedure infinitely many times, Theorem C
says that, in principle, you can know $\Lambda_{\gamma^-}f$.

As an application of our method we consider an inverse problem for
elastic material occupying $\Omega\subset\Bbb R^3$. For general
information about the linear theory of elasticity we refer the
reader to Gurtin's beautiful survey paper \cite{G}. We consider
$\Omega$ as an isotropic elastic body with  Lam\'e parameters
$\lambda, \mu$. In what follows, unless otherwise stated,
$(\lambda,\mu)$ satisfies
$$\left\{
\begin{array}{l}
\text{$\lambda, \mu$ are essentially bounded
functions on $\Omega$ and satisfy }\\
\exists C>0\,\,\text{$\mu(y)>C$ and $3\lambda(y)+2\mu(y)>C$ for almost all $y\in\Omega$.}
\end{array}
\right.
\tag E
$$

Let $Sym A$ denote the symmetric part of the matrix $A$.
One can easily prove that, for each $\mbox{\boldmath $f$}\in\{H^{1/2}(\partial\Omega)\}^3$
there exists a unique $\mbox{\boldmath $u$}\in \{H^1(\Omega)\}^3$ such that
$$\begin{array}{c}
\displaystyle
{\cal L}_{\lambda,\,\mu}\mbox{\boldmath $u$}\equiv(\sum_{1\le j\le 3}\frac{\partial}{\partial y_j}
\{\lambda(\nabla\cdot\mbox{\boldmath $u$})\delta_{ij}+2\mu(Sym\nabla
\mbox{\boldmath $u$})_{ij}\})=0\quad\text{in}\, \Omega,\\
\\
\displaystyle
\mbox{\boldmath $u$}=\mbox{\boldmath $f$}\quad\text{on}\, \partial\Omega.
\end{array}
$$

\noindent
We define the Dirichlet-to-Neumann map, denoted by $\Lambda_{\lambda,\,\mu}$,
by the formula
$$\displaystyle
<\Lambda_{\lambda,\,\mu}\mbox{\boldmath $f$}, \mbox{\boldmath $g$}>
=\int_{\Omega}\lambda\nabla\cdot\mbox{\boldmath $u$}\nabla\cdot
\mbox{\boldmath $v$}
+2\mu Sym\nabla\mbox{\boldmath $u$}\cdot Sym\nabla\mbox{\boldmath $v$}dy
$$
where $\mbox{\boldmath $v$}
\in\{H^1(\Omega)\}^n$ with
$\mbox{\boldmath $v$}
\vert_{\partial\Omega}=
\mbox{\boldmath $g$}
\in\{H^{1/2}(\partial\Omega)\}^3$.
$\Lambda_{\lambda,\,\mu}
\mbox{\boldmath $f$}$ is the traction on $\partial\Omega$ corresponding to
a displacement field $\mbox{\boldmath $f$}$ on $\partial\Omega$.

Let $\Gamma$ be a given nonempty open subset of $\partial\Omega$.
Set
$$\displaystyle
{\cal D}(\Gamma)^3=\{\mbox{\boldmath $f$}
\in \{H^{1/2}(\partial\Omega)\}^3
\,\vert\, \text{supp}\,\mbox{\boldmath $f$}
\subset\Gamma\}.
$$
We call the map
$$\displaystyle
{\cal D}(\Gamma)^3\ni\mbox{\boldmath $f$}
\longmapsto \Lambda_{\lambda,\,\mu}\mbox{\boldmath $f$}
\vert_{\Gamma}
$$
the localized Dirichlet-to-Neumann map.
\par
We assume  that $\Omega$ contains an unknown inclusion "discontinuously"
imbedded in a reference medium with known Lam\'e parameters $\lambda_0, \mu_0$.
We denote by $D$ the shape of the inclusion and by $\lambda\vert_D, \mu\vert_D$
its Lam\'e parameters.
\par
The problem is to find a reconstruction procedure of $D$
by means of the localized Dirichlet-to-Neumann map.  In this problem
we do not assume that $\lambda\vert_D, \mu\vert_D$ are known.
Let us describe the result.  We assume that $\lambda_0, \mu_0$ satisfy
(E) and that
$$\begin{array}{c}
\text{$\partial\Omega$ is $C^3$,}\\
\text{$\lambda_0\in C^2(\overline\Omega)$ and $\mu_0\in C^3(\overline\Omega)$}.
\end{array}
$$
In the sequel we will assume the following hypotheses and notation to be in force.
$$\left\{
\begin{array}{l}
\text{(i) $D$ is an open subset of $\Omega$
with Lipschitz boundary}\\
\text{and satisfies $\overline D\subset\Omega$;}\\
\text{(ii) $\Omega\setminus\overline D$ is connected;}\\
\text{(iii)  $(\lambda(x), \mu(x))=(\lambda_0(x), \mu_0(x))$
for almost all $x\in \Omega\setminus\overline D$;}\\
\text{(iv) $\forall a\in\partial D\,\exists(\alpha, \beta)\forall\epsilon>0\exists\delta>0$}\\
\text{$\vert\lambda(x)-\alpha\vert+\vert\mu(x)-\beta\vert<\epsilon$\,\,
for almost all $x\in B(a,\delta)\cap D$.}\\
\text{Since such $(\alpha, \beta)$ is unique and depends on $(\lambda\vert_D, \mu\vert_D)$
and $a$, we write}\\
\text{$(\lambda_D(a), \mu_D(a))=(\alpha, \beta)$;}\\
\text{(v) $(\lambda_D(a), \mu_D(a))\not=(\lambda_0(a), \mu_0(a))$ for all $a\in\partial D$.}
\end{array}
\right.
\tag {EI}
$$

\noindent
In Section 5 we prove

\proclaim{\noindent Theorem D.}
Assume the conditions of (EI).
Then, there exists a set ${\cal D}(\lambda_0,\mu_0,\Gamma)$
contained in ${\cal D}(\Gamma)^3$ such that
$\partial D$ can be reconstructed from the set
$$
\{<(\Lambda_{\lambda,\,\mu}-\Lambda_{\lambda_0,\,\mu_0})
\mbox{\boldmath $f$},
\mbox{\boldmath $f$}
>
\,\vert\,\mbox{\boldmath $f$}
\in{\cal D}(\lambda_0,\mu_0,\Gamma)\}.
$$
\endproclaim

{\it\noindent Remark.}  Since $\text{supp}\,\mbox{\boldmath $f$}
\subset\Gamma$, we use
only the restriction of $\Lambda_{\lambda, \mu}\mbox{\boldmath $f$}$ to
$\Gamma$ for recovering $D$; we do not require any regularity on
$\lambda, \mu$ inside the inclusion; for the concrete description
of ${\cal D}(\lambda_0, \mu_0, \Gamma)$ see Remark under Proposition
5.6.

The condition (v) of (EI) makes the problem difficult.
Such type of condition never appeared in inverse conductivity problem.
To overcome the difficulty we fully make use  of the speciality of isotropic elastic body.

The second problem is to find a reconstruction formula of
$(\lambda\vert_D, \mu\vert_D)$ by means of the localized
Dirichlet-to-Neumann map provided $D$ is known.  To do it a result
similar to Theorem C shall be useful.

{\bf\noindent Definition(Dirichlet-to-Neumann Map inside Inclusion).}
Consider $D$ and $\lambda, \mu$ satisfying (i), (ii) of (EI) and (E) respectively.
We define the Dirichlet-to-Neumann map inside inclusion
$$\displaystyle
\Lambda_{\lambda^-,\,\mu^-}
:\{H^{1/2}(\partial D)\}^3\longrightarrow \{\{H^{1/2}(\partial D)\}^3\}^*
$$
by the formula
$$\displaystyle
<\Lambda_{\lambda^-,\,\mu^-}\mbox{\boldmath $f$},
\mbox{\boldmath $g$}
>=\int_D\lambda\nabla\cdot
\mbox{\boldmath $u$}
\nabla\cdot
\mbox{\boldmath $v$}
+2\mu Sym\nabla\mbox{\boldmath $u$}
\cdot Sym\nabla\mbox{\boldmath $v$}
dx,
$$
where $\mbox{\boldmath $g$}
\in \{H^{1/2}(\partial D)\}^3$, $\mbox{\boldmath $v$}$ is any $\{H^1(D)\}^3$ function
with $\mbox{\boldmath $v$}
\vert_{\partial D}=\mbox{\boldmath $g$}$, $\mbox{\boldmath $u$}$ is in $\{H^1(D)\}^3$ and
the weak solution of the Dirichlet problem
$$\begin{array}{c}
\displaystyle
{\cal L}_{\lambda,\,\mu}\mbox{\boldmath $u$}
=0\,\,\text{in}\, D,\\
\\
\displaystyle
\mbox{\boldmath $u$}
\vert_{\partial D}=\mbox{\boldmath $f$}.
\end{array}
$$

The proof of Theorem E stated below proceeds along the same lines
with that of Theorem C and we omit it.

\proclaim{\noindent Theorem E.}  Let
$\Gamma$ be a given nonempty open subset of $\partial\Omega$.
Assume that (i)$\sim$(iii) of (EI) hold. There exists a set
${\cal D}(\lambda_0, \mu_0, \Gamma, D)$ contained in ${\cal D}(\Gamma)^3$
such that $\Lambda_{\lambda^{-},\,\mu^{-}}$ can be calculated from
the set
$$\displaystyle
\{<(\Lambda_{\lambda,\,\mu}-\Lambda_{\lambda_0,\,\mu_0})\mbox{\boldmath $f$}
,\mbox{\boldmath $f$}
>\,\vert
\, \mbox{\boldmath $f$}
\in{\cal D}(\lambda_0, \mu_0, \Gamma,D)\}.
$$
\endproclaim

{\it\noindent Remark.}  In this theorem, it is assumed that
$\lambda_0, \mu_0$ and $D$ are known;  ${\cal D}(\gamma_0, \mu_0,\Gamma, D)$
is independent of $\Lambda_{\lambda,\,\mu}$ and depends
on $\Omega, \Gamma$ and $D$;
notice that we do not assume any regularity on $\lambda, \mu$ inside the inclusion.

By the way, Nakamura-Uhlmann \cite{NU1} proved

\proclaim{\noindent Theorem F.}  Let $(\lambda_i, \mu_i), i=1, 2, $
be  two pairs of Lam\'e parameters.
Assume that
$$
\text{$\lambda_i, \mu_i$ are smooth on $\overline\Omega$,\,$i=1, 2$,}
$$
and
$$
\text{$\partial\Omega$ is smooth.}
$$
If
$$
\Lambda_{\lambda_1,\,\mu_1}=\Lambda_{\lambda_2,\,\mu_2},
$$
then
$$
(\lambda_1, \mu_1)=(\lambda_2, \mu_2).
$$
\endproclaim

In contrast to what Theorem A tells Theorem F is a uniqueness
theorem and does not contain  any reconstruction procedure. So far
such procedure is not known. This affects Theorem G mentioned
below.  We refer the reader to their survey paper \cite{NU2} for
more precise information about related results.

We consider two arbitrary triples $(D_1, \lambda_1, \mu_1),
(D_2, \lambda_2, \mu_2)$ both satisfying (EI).
From Theorems D, E and F we immediately get

\proclaim{\noindent Theorem G.}
Assume that for $i=1,2$
$$
\text{$\partial D_i$ is smooth}
$$
and
$$
\text{$\lambda_i\vert_{D_i}, \mu_i\vert_{D_i}$ are smooth on $\overline D_i$.}
$$
There exists a set ${\cal D}(\lambda_0, \mu_0, \Gamma)$ contained
in ${\cal D}(\Gamma)^3$ such that if
$$\displaystyle
\Lambda_{\lambda_1,\,\mu_1}\mbox{\boldmath $f$}
\vert_{\Gamma}
=\Lambda_{\lambda_2,\,\mu_2}\mbox{\boldmath $f$}
\vert_{\Gamma}
$$
for all $\mbox{\boldmath $f$}
\in{\cal D}(\lambda_0, \mu_0, \Gamma)$, then
$$
\text{$D_1=D_2$ and $(\lambda_1\vert_{D_1}, \mu_1\vert_{D_1})
=(\lambda_2\vert_{D_2}, \mu_2\vert_{D_2})$.}
$$
\endproclaim

This is a uniqueness theorem.
Nakamura-Uhlmann's result  cannot cover Theorem G without using Theorem E
since their method heavily relies on the regularity
of $\lambda, \mu$.

\section{The probe method and proof of Theorem B}

Let us recall some definitions in \cite{IK2}.
We insert a ``needle" into $\Omega$ defined below.

{\bf\noindent Definition(Needle).}
We call a continuous map $c:[0, 1]\longrightarrow \overline\Omega$ satisfying
(i) and (ii) a needle:
$$\begin{array}{c}
(i)\quad c(0),\quad c(1) \in \partial\Omega\\
(ii)\quad c(t)\in\Omega\quad (0<t<1).
\end{array}
$$

{\bf\noindent Definition(Impact parameter).}
It is easy to verify that if a needle $c$ touches a point on
$\partial D$, there exists a unique $t(c;D)\in]0, 1[$ such that if
$0<t<t(c;D)$, $c(t)\in\Omega\setminus\overline D$ and
$c(t(c;D))\in \partial D$.  We set $t(c;D)=1$ if $c$ does not
touch any point on $\partial D$.  We call $t(c;D)$ the impact
parameter of $c$ with respect to $D$.  Notice that $t(c;D)$ has
the form
$$\displaystyle
t(c;D)=\sup\{s\in]0, 1[\,\vert\,\forall t\in]0, s[\, c(t)\in\Omega\setminus\overline D\}.
\tag {2.1}
$$

By virtue of (ii) of (CI) we get
$$\displaystyle
\partial D=\{c(t(c;D))\,\vert\, \text{$c$ is a needle and satisfies $t(c;D)<1$}\}.
\tag {2.2}
$$
Our purpose is to construct the analytic version of the impact parameter by means
of the localized Dirichlet-to-Neumann map and show the coincidence of both
impact parameters.  Then from (2.2) we get $\partial D$.

Let us explain its construction.
First we have to prove

\proclaim{\noindent Proposition 2.1.}
For each $x\in\Omega$ let $G_0(\cdot;x)$ denote the standard
fundamental solution for the operator $\nabla\cdot\gamma_0(x)\nabla$.
Then, there exists a family $(G_x^0(\cdot))_{x\in\Omega}$
in $\cap_{1\le p<\frac{n}{n-2}}L^p(\Omega)$ such that
$$
\text{$(G_x^0(\cdot)-G_0(\cdot-x;x))_{x\in\Omega}$ is bounded
in $H^1(\Omega)$,}
$$
$$
\displaystyle
\nabla\cdot\gamma_0\nabla G_x^0(\cdot)+\delta(\cdot-x)=0\,\,\text{in}\,\Omega.
$$
\endproclaim

For the proof see Section 6.
Second we prepare

\proclaim{\noindent Proposition 2.2(Selection of boundary data).}  Let $c$ be a needle.
Then, for any $t\in]0, 1[$, there exists a sequence
$\{f_n(\cdot ; c(t))\}$ of ${\cal D}(\Gamma)$ such that the weak solution $v_n$ of
$$\begin{array}{c}
\displaystyle
\nabla\cdot\gamma_0\nabla v=0\quad\text{in}\,\Omega,\\
\\
\displaystyle
v\vert_{\partial\Omega}=f_n(\cdot ;c(t))
\end{array}
$$
converges to $G_{c(t)}^0(\cdot)$ in $H^1_{\text{loc}}(\Omega\setminus\{c(t')\vert
0<t'\le t \})$ as $n\longrightarrow\infty$.
\endproclaim

{\it\noindent Proof.}
By virtue of the connectedness of $\partial\Omega$,
we can take a sequence $\{O_n\}$ of relatively compact open subsets of
$\Omega\setminus\{c(t')\vert 0<t'\le t \}$ such that
$$\begin{array}{c}
\overline O_n\subset O_{n+1},\\
\\
\displaystyle
\Omega\setminus\{c(t')\vert 0<t'\le t \}=\cup_n O_n,\\
\\
\displaystyle
\text{$\Omega\setminus\overline O_n$ is connected.}
\end{array}
$$
Then, for each $O_n$,
apply the Runge approximation property in Appendix.

\noindent
$\Box$

In \cite{IK1} the author proved

\proclaim{\noindent Proposition 2.3.}
Let $\gamma_j, j =1, 2, $ be two conductivities.
Let $v_j\in H^1(\Omega)$ denote the weak solution of the Dirichlet problem
$$\begin{array}{c}
\displaystyle
\nabla\cdot\gamma_j\nabla v=0\,\,\text{in}\,\Omega,\\
\\
\displaystyle
v\vert_{\partial\Omega}=f.
\end{array}
$$
It holds that
$$\begin{array}{c}
\displaystyle
\int_{\Omega}\{\gamma_1^{-1}-\gamma_2^{-1}\}\gamma_1\nabla v_1\cdot\gamma_1
\nabla v_1dx
\le <(\Lambda_{\gamma_2}-\Lambda_{\gamma_1})f\,,f>,\\
\\
\displaystyle
<(\Lambda_{\gamma_2}-\Lambda_{\gamma_1})f\,,f>
\le\int_{\Omega}(\gamma_2-\gamma_1)\nabla v_1\cdot\nabla v_1 dx,\\
\\
\displaystyle
\int_{\Omega}\{\gamma_2^{-1}-\gamma_1^{-1}\}\gamma_2\nabla v_2\cdot\gamma_2\nabla v_2dx
\le <(\Lambda_{\gamma_1}-\Lambda_{\gamma_2})f\,,f>,\\
\\
\displaystyle
<(\Lambda_{\gamma_1}-\Lambda_{\gamma_2})f,\,f>
\le\int_{\Omega}(\gamma_1-\gamma_2)\nabla v_2\cdot\nabla v_2 dx.
\end{array}
\tag {2.3}
$$
\endproclaim

{\it\noindent Remark.}
These however are not the best possible estimates.
In fact, therein the author proved also
$$\begin{array}{c}
\displaystyle
\int_{\Omega}\{\gamma_1^{-1}-\gamma_2^{-1}\}\gamma_1\nabla v_1
\cdot\gamma_1\nabla v_1 dx
\le \frac{<\Lambda_{\gamma_1}f, f>}{<\Lambda_{\gamma_2}f, f>}
<(\Lambda_{\gamma_2}-\Lambda_{\gamma_1})f, f>,\\
\\
\displaystyle
\int_{\Omega}\{\gamma_2^{-1}-\gamma_1^{-1}\}\gamma_2\nabla v_2
\cdot\gamma_2\nabla v_2 dx
\le \frac{<\Lambda_{\gamma_2}f, f>}{<\Lambda_{\gamma_1}f, f>}
<(\Lambda_{\gamma_1}-\Lambda_{\gamma_2})f, f>.
\end{array}
$$
The proof is not so trivial.  Notice that if $\gamma_1$ and $ \gamma_2$ are constant,
the inequalities indicated above become the equalities.

These inequalities can be derived just from the definition of the
Dirichlet-to-Neumann map like as the Schwartz inequality in
Hilbert space theory. For complete details, we refer the reader to
\cite{IK1} or \cite{IK2}.

{\bf\noindent Definition($T(c)$ and $I(t,c)$).}
Denote the set of all $s\in]0, 1[$ such that
$$\displaystyle
I(t,c)\equiv\lim_{n\rightarrow\infty}
 <(\Lambda_{\gamma}-\Lambda_{\gamma_0})f_n(\cdot;c(t))\,,f_n(\cdot;c(t))>
$$
exists for $0<t<s$ and $\sup_{0<t<s}\vert I(t,c)\vert<\infty$
by $T(c)$.

{\it\noindent Remark.}
The definition of $T(c)$ is slightly different from that of \cite{IK3}.

A combination of Proposition 2.2 and (2.3) gives

\proclaim{\noindent Proposition 2.4.}
Let $c$ be a needle.
If $0<t<t(c;D)$, $I(t,c)$
exists and it holds that
$$\begin{array}{c}
\displaystyle
\int_D\{\gamma_0^{-1}-\gamma^{-1}\}\gamma_0\nabla G_{c(t)}^0\cdot
\gamma_0\nabla G_{c(t)}^0dx\\
\\
\displaystyle
\le I(t,c) \le\int_D (\gamma-\gamma_0)\nabla G_{c(t)}^0\cdot\nabla G_{c(t)}^0dx.
\end{array}
\tag {2.4}
$$
\endproclaim

{\it\noindent Proof.}  Let $u_n\in H^1(\Omega)$ be the weak solution of
$$\begin{array}{c}
\displaystyle
\nabla\cdot\gamma\nabla u_n=0\,\,\text{in}\, \Omega,\\
\\
\displaystyle
u_n\vert_{\partial\Omega}=f_n(\cdot;c(t)).
\end{array}
$$
Since $w_n=u_n-v_n\in H^1_0(\Omega)$ satisfies
$$
\nabla\cdot\gamma\nabla w_n=-\nabla\cdot\chi_D(\gamma-\gamma_0)\nabla v_n\,\,\text{in}\, \Omega
$$
and
$$\displaystyle
v_n\longrightarrow G^0_{c(t)}(\cdot)\,\,\text{in}\, H^1(D)\,\,\text{for}\, 0<t<t(c;D),
$$
we know that $w_n$ converges  in $H^1(\Omega)$ to the weak solution of
$$\begin{array}{c}
\displaystyle
\nabla\cdot\gamma\nabla w=-\nabla\cdot\chi_D(\gamma-\gamma_0)\nabla G^0_{c(t)}\,\,
\text{in}\, \Omega,\\
\\
\displaystyle
w\vert_{\partial\Omega}=0.
\end{array}
$$
Then from Alessandini's identity we get
$$\begin{array}{c}
\displaystyle
I(t,c)=\lim_{n\longrightarrow\infty}\int_{D}(\gamma-\gamma_0)\nabla u_n\cdot\nabla v_n dx\\
\\
\displaystyle
=\int_D(\gamma-\gamma_0)\nabla(G^0_{c(t)}+w)\cdot\nabla G^0_{c(t)} dx.
\end{array}
\tag {2.5}
$$
(2.4) is clearly valid.

\noindent
$\Box$

Theorem B is a consequence of

\proclaim{\noindent Proposition 2.5(Reconstruction of impact parameter).}
\newline{For any needle $c$} it holds that
$$\displaystyle
T(c)=]0, t(c;D)[
$$
and thus the formula
$$\displaystyle
t(c;D)=\sup T(c),
\tag {2.6}
$$
is valid.  $\partial D$ has the form
$$\displaystyle
\partial D=\{c(\sup T(c))\,\vert\, \text{$c$ is a needle and satisfies
$\sup T(c)<1$}\}.
\tag {2.7}
$$
\endproclaim

{\it Remark.}  ${\cal D}(\gamma_0, \Gamma)$ in Theorem B is
$$\displaystyle
\{f_n(\cdot;c(t))\,,n=1,\cdots\,\vert\,\text{$c$ is a needle and $0<t<1$}\}.
$$

The proof of Proposition 2.5 proceeds along the almost same lines with that of Theorem A in \cite{IK3}
and so we describe its outline.
First from Propositions 2.1 and 2.4 we get $]0, t(c;D)[\subset T(c)$.
Since $T(c)\subset ]0, 1[$, if $t(c;D)=1$, we know that $T(c)=]0, 1[=]0, t(c;D)[$.
So the problem is the case where $t(c;D)<1$.  Set $a=c(t(c;D))(\in\partial D)$.
We assume that $]0, t(c;D)[$ does not coincide with $T(c)$.  Then there exists $s\in T(c)$
such that $t(c;D)\le s$.  Then $\sup _{0<t<t(c;D)}\vert I(t,c)\vert<C$ for a constant.
But a combination of (iv) of (CI), Propositions 2.1 and 2.4 gives
$$\displaystyle
\lim_{t\uparrow t(c;D)}\vert I(t,c)\vert=\infty.
$$
This is a contradiction.

Proposition 2.5 tells us that
using the localized Dirichlet-to-Neumann map, we can predict whether
you encounter a point on the boundary of unknown inclusion when
you go down along the given needle from $\partial\Omega$.

\section{Proof of Theorem C}

For each $f\in H^{-1/2}(\partial D)$ define
$$\displaystyle
T_f(\varphi)=<f,\,\varphi\vert_{\partial D}>, \varphi\in H^1_0(\Omega).
$$
From the trace theorem we know $T_f\in H^{-1}(\Omega)\equiv (H^1_0(\Omega))^*$
and thus there exists a unique weak solution $u_f$ in $H^1_0(\Omega)$ of
$$\begin{array}{c}
\displaystyle
\nabla\cdot\gamma\nabla u_f=-T_f\,\,\text{in}\, \Omega,\\
\\
\displaystyle
u_f\vert_{\partial\Omega}=0.
\end{array}
\tag {3.1}
$$
Set
$$\displaystyle
Gf=u_f\vert_{\partial D}.
$$
It is easy to see that $G$ is a bounded linear operator
from $H^{-1/2}(\partial D)$ to $H^{1/2}(\partial D)$.

{\bf\noindent  Definition(Dirichlet-to-Neumann map outside inclusion).}
We define the Dirichlet-to-Neumann map outside inclusion
$$\displaystyle
\Lambda_{\gamma^+}:H^{1/2}(\partial D)\longrightarrow H^{-1/2}(\partial D)
$$
by the formula
$$\displaystyle
<\Lambda_{\gamma^+}f, g>=-\int_{\Omega\setminus\overline D}
\gamma\nabla u\cdot\nabla\varphi dx,
$$
where $g\in H^{1/2}(\partial D)$, $\varphi$ is any
$H^1(\Omega\setminus\overline D)$ function
with $\varphi\vert_{\partial D}=g$ and $\varphi\vert_{\partial\Omega}=0$,
$u$ is in $H^1(\Omega\setminus\overline D)$ and
the weak solution of the Dirichlet problem
$$\begin{array}{c}
\nabla\cdot\gamma\nabla u=0\,\,\text{in}\, \Omega\setminus\overline D,\\
\\
\displaystyle
u\vert_{\partial \Omega}=0,\\
\\
\displaystyle
u\vert_{\partial D}=f.
\end{array}
$$

We start with

\proclaim{\noindent Proposition 3.1.}
\newline{(i)  $\Lambda_{\gamma^{-}}-\Lambda_{\gamma^{+}}$ is injective;}
\newline{(ii)  the formula}
$$
(\Lambda_{\gamma^-}-\Lambda_{\gamma^+})Gf=f, \forall f\in H^{-1/2}(\partial D)
$$
is valid.
\endproclaim

From (i) and (ii) of Proposition 3.1 we can conclude that
$\Lambda_{\gamma^{-}}-\Lambda_{\gamma^{+}}$ is bijective
and thus we know $G=(\Lambda_{\gamma^-}-\Lambda_{\gamma^+})^{-1}$.
Therefore $G$ is bijective, too and hence the formula
$$\displaystyle
\Lambda_{\gamma^-}-\Lambda_{\gamma^+}=G^{-1},
\tag {3.2}
$$
is valid.  Nachman ((6.15) in \cite{N2}) proved the corresponding fact
in the case where $\Omega\subset\Bbb R^2$ and $\gamma\in W^{2,p}(\Omega)$
with $p>1$.  Our proof is quite elementary than that
of Nachman and notice that we do not assume any regularity of $\gamma$ on $\Omega$.
This is because of the weak formulation of $G$.

{\it\noindent Proof of Proposition 3.1.}
(i) is well known, and is proved here only for the convenience of the reader.
Assume that $g\in H^{1/2}(\partial D)$
satisfies $(\Lambda_{\gamma^-}-\Lambda_{\gamma^+})g=0$.
Let $u^-\in H^1(D)$ be the weak solution of
$$\begin{array}{c}
\displaystyle
\nabla\cdot\gamma\nabla u^-=0\,\,\text{in}\, D,\\
\\
\displaystyle
u^-\vert_{\partial D}=g.
\end{array}
$$
Let $u^+\in H^1(\Omega\setminus\overline D)$ be the weak solution of
$$\begin{array}{c}
\displaystyle
\nabla\cdot\gamma\nabla u^+=0\,\,\text{in}\,\Omega\setminus\overline D,\\
\\
\displaystyle
u^+\vert_{\partial D}=g,\\
\\
\displaystyle
u\vert_{\partial\Omega}=0.
\end{array}
$$
Set
$$
u=\left\{
\begin{array}{lr}
u^- & \text{in $D$}\\
\\
\displaystyle
u^+ & \text{ in $\Omega\setminus\overline D$}.
\end{array}
\right.
$$
Then $u\in H^1_0(\Omega)$ and from the assumption on $g$ we know that
$u$ is the weak solution of
$$\begin{array}{c}
\displaystyle
\nabla\cdot\gamma\nabla u=0\,\,\text{in}\, \Omega,\\
\\
\displaystyle
u\vert_{\partial\Omega}=0.
\end{array}
$$
Therefore $u=0$ and thus $g=0$.

Let us give the proof of (ii).  Let $\varphi\in H^1_0(\Omega)$.
From (3.1) we get
$$\begin{array}{c}
\displaystyle
<f,\,\varphi\vert_{\partial D}>=T_f(\varphi)\\
\\
\displaystyle
=\int_{\Omega}\gamma\nabla u_f\cdot\nabla\varphi dx\\
\\
\displaystyle
=\int_{\Omega\setminus\overline D}\gamma\nabla u_f\cdot\nabla\varphi dx
+\int_{D}\gamma\nabla u_f\cdot\nabla\varphi dx\\
\\
\displaystyle
=-<\Lambda_{\gamma^+}(u_f\vert_{\partial D})\,,\varphi\vert_{\partial D}>
+<\Lambda_{\gamma^-}(u_f\vert_{\partial D})\,,\varphi\vert_{\partial D}>\\
\\
\displaystyle
=<(\Lambda_{\gamma^-}-\Lambda_{\gamma^+})Gf\,,\varphi\vert_{\partial D}>.
\end{array}
$$
From this identity we get (ii) since the map
$$\displaystyle
H^1_0(\Omega)\ni\varphi\longmapsto\varphi\vert_{\partial D}\in H^{1/2}(\partial D)
$$
is surjective.

\noindent
$\Box$

In the remainder of this section we make use of (i)$\sim$(iii) of
(CI).  Let $F$ be a given
element of $H^{-1}(\Omega)$ with
$$\displaystyle
\text{supp}\,F\subset\Omega\setminus\overline D.
$$
Let $u_0\in H^1_0(\Omega)$ be the weak solution of
$$\begin{array}{c}
\displaystyle
\nabla\cdot\gamma_0\nabla u_0=-F\,\,\text{in}\, \Omega,\\
\\
\displaystyle
u_0\vert_{\partial\Omega}=0
\end{array}
\tag {3.3}
$$
Set
$$\displaystyle
\mbox{\boldmath $G$}_0F=u_0.
$$
Let $w\in H^1_0(\Omega)$ be the weak solution of
$$\begin{array}{c}
\displaystyle
\nabla\cdot\gamma\nabla w=-\nabla\cdot\chi_D(\gamma-\gamma_0)\nabla u_0
\,\,\text{in}\, \Omega,\\
\\
\displaystyle
w\vert_{\partial\Omega}=0.
\end{array}
\tag {3.4}
$$
Set
$$\displaystyle
\mbox{\boldmath $G$}
F=\mbox{\boldmath $G$}_0F+w.
$$

\noindent
We prove the one of two crucial identities, needed for calculating
$G$.

\proclaim{\noindent Proposition 3.2.}  The formula
$$\displaystyle
F(u_f)=<f, \mbox{\boldmath $G$}
F\vert_{\partial D}>,
\tag {3.5}
$$
is valid.
\endproclaim

{\it\noindent Proof.}  A combination of (3.1) and (3.4) yields
$$\begin{array}{c}
\displaystyle
\int_D(\gamma-\gamma_0)\nabla u_f\cdot\nabla u_0 dx
=-\int_{\Omega}\gamma\nabla w\cdot\nabla u_f dx\\
\\
\displaystyle
=-T_f(w)\\
\\
\displaystyle
=-<f, w\vert_{\partial D}>.
\end{array}
\tag {3.6}
$$
On the other hand, from (3.1) and (3.3) we get
$$\begin{array}{c}
<f, u_0\vert_{\partial D}>
=T_f(u_0)\\
\\
\displaystyle
=\int_{\Omega}\gamma\nabla u_f\cdot\nabla u_0 dx\\
\\
\displaystyle
=\int_{\Omega}(\gamma-\gamma_0)\nabla u_f\cdot\nabla u_0 dx
+\int_{\Omega}\gamma_0\nabla u_f\cdot\nabla u_0 dx\\
\\
\displaystyle
=\int_{D}(\gamma-\gamma_0)\nabla u_f\cdot\nabla u_0 dx
+F(u_f).
\end{array}
$$
Combining this with (3.6), we get (3.5).

\noindent
$\Box$

Let $H$ be a given
element of $H^{-1}(\Omega)$ with
$$
\text{supp}\,H\subset\Omega\setminus\overline D.
$$
Let $v_0\in H^1_0(\Omega)$ be the weak solution of
$$\begin{array}{c}
\displaystyle
\nabla\cdot\gamma_0\nabla v_0=-H\,\,\text{in}\, \Omega,\\
\\
\displaystyle
v_0\vert_{\partial\Omega}=0.
\end{array}
\tag {3.7}
$$
Notice that both $u_0$ and $v_0$ satisfies $\nabla\cdot\gamma_0\nabla u=0$
in an open neighbourhood of $\overline D$.  Therefore from the Runge approximation
property proved in Appendix we get two sequences $\{u_n\}$, $\{v_n\}$
of $H^1(\Omega)$ functions with
$$
\text{supp}\,(u_n\vert_{\partial\Omega}), \text{supp}\,(v_n\vert_{\partial\Omega})\subset\Gamma
$$
such that
$$\begin{array}{c}
\displaystyle
\nabla\cdot\gamma_0\nabla u_n=0\,\,\text{in}\,\Omega,\\
\\
\displaystyle u_n\longrightarrow u_0\,\,\text{in}\, H^1(D),\\
\\
\displaystyle
\nabla\cdot\gamma_0\nabla v_n=0\,\,\text{in}\,\Omega,\\
\\
\displaystyle
v_n\longrightarrow v_0\,\,\text{in}\,H^1(D).
\end{array}
$$

\proclaim{\noindent Proposition 3.3.}  The formula
$$\displaystyle
-H(\mbox{\boldmath $G$}
F-\mbox{\boldmath $G$}_0F)
=
\displaystyle
\frac{1}{4}\lim_{n\longrightarrow\infty}
\{<(\Lambda_{\gamma}-\Lambda_{\gamma_0})f_n^+,
f_n^+> -<(\Lambda_{\gamma}-\Lambda_{\gamma_0})f_n^-, f_n^->\},
\tag {3.8}
$$
is valid where
$$\displaystyle
f_n^+=(u_n+v_n)\vert_{\partial\Omega},\, f_n^-=(u_n-v_n)\vert_{\partial\Omega}.
$$
\endproclaim

{\it\noindent Remark.}  ${\cal D}(\gamma_0, \Gamma, D)$ in Theorem C is
$$\displaystyle
\{f_n^+, f_n^-, n=1,\cdots\,\vert\,\text{$F, H\in H^{-1}(\Omega)$ with $\text{supp}\,F,
\text{supp}\,H\subset\Omega\setminus\overline D$}\}.
$$

{\it\noindent Proof of Proposition 3.3.}
Let $a_n\in H^1(\Omega)$ be the weak solution of
$$\begin{array}{c}
\displaystyle
\nabla\cdot\gamma\nabla a_n=0\,\,\text{in}\, \Omega,\\
\\
\displaystyle
a_n\vert_{\partial\Omega}=u_n\vert_{\partial\Omega}.
\end{array}
$$
Set
$$\displaystyle
w_n=a_n-u_n\in H^1_0(\Omega).
$$
$w_n$ solves
$$
\displaystyle
\nabla\cdot\gamma\nabla w_n=-\nabla\cdot\chi_D(\gamma-\gamma_0)
\nabla u_n\,\,\text{in}\, \Omega.
$$
It is easy to see that $w_n\longrightarrow w$ in $H^1(\Omega)$.
Thus from Alessandrini's identity  and the symmetry of the Dirichlet-to-Neumann map we get
$$\begin{array}{c}
\displaystyle
\frac{1}{4}\{<(\Lambda_{\gamma}-\Lambda_{\gamma_0})f_n^+, f_n^+>-
<(\Lambda_{\gamma}-\Lambda_{\gamma_0})f_n^-, f_n^->\}\\
\\
\displaystyle
=<(\Lambda_{\gamma}-\Lambda_{\gamma_0})u_n\vert_{\partial\Omega},
v_n\vert_{\partial\Omega}>\\
\\
\displaystyle
=\int_{\Omega}(\gamma-\gamma_0)\nabla a_n\cdot\nabla v_n dx\\
\\
\displaystyle
=\int_{D}(\gamma-\gamma_0)\nabla a_n\cdot\nabla v_n dx\\
\\
\displaystyle
=\int_{D}(\gamma-\gamma_0)\nabla u_n\cdot\nabla v_n dx
+\int_{D}(\gamma-\gamma_0)\nabla w_n\cdot\nabla v_n dx\\
\\
\displaystyle
\longrightarrow\int_{D}(\gamma-\gamma_0)\nabla u_0\cdot\nabla v_0 dx
+\int_{D}(\gamma-\gamma_0)\nabla w\cdot\nabla v_0 dx.
\end{array}
\tag {3.9}
$$
From (3.4) and (3.7) we get
$$\begin{array}{c}
H(w)=\int_{\Omega}\gamma_0\nabla v_0\cdot\nabla w dx\\
\\
\displaystyle
=\int_{\Omega}\{\gamma-\chi_D(\gamma-\gamma_0)\}\nabla v_0\cdot\nabla wdx\\
\\
\displaystyle
=\int_{\Omega}\gamma\nabla v_0\cdot\nabla w dx-\int_{D}(\gamma-\gamma_0)
\nabla v_0\cdot\nabla w dx\\
\\
\displaystyle
=-\int_{D}(\gamma-\gamma_0)\nabla u_0\cdot\nabla v_0 dx
-\int_{D}(\gamma-\gamma_0)\nabla v_0\cdot\nabla w dx.
\end{array}
$$
Combining this with (3.9), we get (3.8).

\noindent
$\Box$

{\it\noindent Proof of Theorem C.}
First from (3.8) we get $\mbox{\boldmath $G$}
F-\mbox{\boldmath $G$}_0F$ in $\Omega\setminus\overline D$
and thus its trace on $\partial D$.  From (3.5) we get $u_f$ in $\Omega\setminus\overline D$ and thus its trace on $\partial D$, that is $Gf$.  Therefore we get $G$ and from (3.2) $\Lambda_{\gamma^-}$, too.

\noindent
$\Box$

\section{Summary of reconstruction procedure}

\proclaim{Reconstruction of the shape of inclusion.}

(1)  First construct a family $(G_x^0(\cdot))_{x\in\Omega}$ each of which is a special solution
of \newline{$\nabla\cdot\gamma_0\nabla u=0$} in $\Omega\setminus\{x\}$.

(2) For each needle $c$ and $t\in]0, 1[$, take a sequence $\{f_n(\cdot;c(t))\}$ of functions
on $\partial\Omega$ with $\text{supp}\,f_n\subset\Gamma$ in such a way that the solution $v_n$ of
$$\begin{array}{c}
\displaystyle
\nabla\cdot\gamma_0\nabla v=0\,\,\text{in}\, \Omega\\
\\
\displaystyle
v\vert_{\partial\Omega}=f_n
\end{array}
$$
converges to $G_{c(t)}^0$ in $H^1_{\text{loc}}(\Omega\setminus\{c(t')\vert
0<t'\le t\})$ as $n\longrightarrow\infty$.

(3) Calculate $T(c)$.

(4)  Use formula $t(c;D)=\sup T(c)$ to recover $t(c;D)$ from $T(c)$.

(5) Use formula $\partial D=\{c(t(c;D))\vert \text{$c$ is a needle and $t(c;D)<1$}\}$ to recover $\partial D$ from $t(c;D)$.

\endproclaim

\proclaim{\noindent Reconstruction of the Dirichlet-to-Neumann map inside inclusion.}

(1)  Give $F, H$ being elements of $H^{-1}(\Omega)$ with $\text{supp}\,F,\,\text{supp}\,H\subset\Omega\setminus\overline D$.

(2)  Construct a sequence $\{u_n\}$ of functions in $H^1(\Omega)$ such that
$$\begin{array}{c}
\displaystyle
\nabla\cdot\gamma_0\nabla u_n=0\,\,\text{in}\, \Omega,\\
\\
\displaystyle
\text{supp}\,(u_n\vert_{\partial\Omega})\subset\Gamma,\\
\\
\displaystyle
u_n\longrightarrow \mbox{\boldmath $G$}_0F\,\,\text{in}\, H^1(D).
\end{array}
$$

(3)
 Construct a sequence $\{v_n\}$ of functions in $H^1(\Omega)$ such that
$$\begin{array}{c}
\nabla\cdot\gamma_0\nabla v_n=0\,\,\text{in}\, \Omega,\\
\\
\displaystyle
\text{supp}\,(v_n\vert_{\partial\Omega})\subset\Gamma,\\
\\
\displaystyle
v_n\longrightarrow \mbox{\boldmath $G$}
_0H\,\,\text{in}\, H^1(D).
\end{array}
$$

(4)  Calculate
$$\displaystyle
f_n^+=(u_n+v_n)\vert_{\partial\Omega}\,,f_n^-=(u_n-v_n)\vert_{\partial\Omega}.
$$

(5) Use formula
$$\displaystyle
-H(\mbox{\boldmath $G$}
F-\mbox{\boldmath $G$}_0F)
=\frac{1}{4}\lim_{n\longrightarrow\infty}
\{<(\Lambda_{\gamma}-\Lambda_{\gamma_0})f_n^+,
f_n^+> -<(\Lambda_{\gamma}-\Lambda_{\gamma_0})f_n^-, f_n^->\}
$$
to recover $H(\mbox{\boldmath $G$}
F-\mbox{\boldmath $G$}_0F)$ from the set
$$\displaystyle
\{<(\Lambda_{\gamma}-\Lambda_{\gamma_0})f_n^+, f_n^+>,
<(\Lambda_{\gamma}-\Lambda_{\gamma_0})f_n^-, f_n^->\}.
$$

(6) Calculate $\mbox{\boldmath $G$}
F-\mbox{\boldmath $G$}_0F$ in $\Omega\setminus\overline D$
from the set
$$
\{H(\mbox{\boldmath $G$}
F-\mbox{\boldmath $G$}_0F)\,\vert\,
\text{$H\in H^{-1}(\Omega)$ with $\text{supp}\,H\subset\Omega\setminus\overline D$}\}.
$$

(7)  Calculate $\mbox{\boldmath $G$}
F\vert_{\partial D}$.

(8)  Give $f\in H^{-1/2}(\partial D)$.

(9)  Use formula
$$
F(u_f)=<f, \mbox{\boldmath $G$}
F\vert_{\partial D}>
$$
to recover $F(u_f)$ from $\mbox{\boldmath $G$}
F\vert_{\partial D}$.

(10)  Calculate $u_f$ in $\Omega\setminus\overline D$ from the set
$$
\{F(u_f)\,\vert\,\text{$F\in H^{-1}(\Omega)$
 with $\text{supp}\,F\subset\Omega\setminus\overline D$}\}.
$$

(11) Calculate $Gf=u_f\vert_{\partial D}$.

(12)  Calculate $G$ from the set
$$
\{Gf\,\vert\, f\in H^{-1/2}(\partial D)\}.
$$

(13)  Use formula
$$
\Lambda_{\gamma^-}-\Lambda_{\gamma^+}=G^{-1}
$$
to recover $\Lambda_{\gamma^-}$ from $G$.

\endproclaim

\section{The multiprobe method and proof of Theorem D}

The starting point is

\proclaim{\noindent Proposition 5.1.}  Let $(\lambda_j, \mu_j), j=1, 2, $ be two pairs of Lam\'e parameters.
Let $\mbox{\boldmath $u$}
_j\in\{H^1(\Omega)\}^3$ denote the weak solution of
$$\begin{array}{c}
\displaystyle
{\cal L}_{\lambda_j,\,\mu_j}\mbox{\boldmath $u$}
_j=0\quad\text{in $\Omega$}\\
\\
\displaystyle
\mbox{\boldmath $u$}
_j\vert_{\partial\Omega}=\mbox{\boldmath $f$}.
\end{array}
$$
It holds that
$$\begin{array}{c}
\displaystyle
\int_{\Omega}\frac{3\lambda_2+2\mu_2}{3(3\lambda_1+2\mu_1)}
\{3(\lambda_1-\lambda_2)+2(\mu_1-\mu_2)\}\vert\nabla\cdot\mbox{\boldmath $u$}
_2\vert^2\\
\\
\displaystyle
+\frac{\mu_2}{\mu_1}2(\mu_1-\mu_2)\left\vert Sym\nabla\mbox{\boldmath $u$}
_2-\frac{\nabla\cdot\mbox{\boldmath $u$}
_2}{3}I_3
\right\vert^2 dx\\
\\
\displaystyle
\le <(\Lambda_{\lambda_1,\,\mu_1}-\Lambda_{\lambda_2,\,\mu_2})\mbox{\boldmath $f$}
, \mbox{\boldmath $f$}>,
\end{array}
\tag {5.1}
$$
$$\begin{array}{c}
\displaystyle
<(\Lambda_{\lambda_1,\,\mu_1}-\Lambda_{\lambda_2,\,\mu_2})\mbox{\boldmath $f$}
, \mbox{\boldmath $f$}>\\
\\
\displaystyle
\le\int_{\Omega}\frac{3(\lambda_1-\lambda_2)+2(\mu_1-\mu_2)}{3}
\left\vert\nabla\cdot\mbox{\boldmath $u$}
_2\right\vert^2\\
\\
\displaystyle
+2(\mu_1-\mu_2)\left\vert Sym\nabla\mbox{\boldmath $u$}
_2
-\frac{\nabla\cdot\mbox{\boldmath $u$}
_2}{3}I_3\right\vert^2 dx,
\end{array}
\tag {5.2}
$$
$$\begin{array}{c}
\displaystyle
\int_{\Omega}\frac{3\lambda_1+2\mu_1}{3(3\lambda_2+2\mu_2)}
\{3(\lambda_2-\lambda_1)+2(\mu_2-\mu_1)\}\vert\nabla\cdot\mbox{\boldmath $u$}
_1\vert^2\\
\\
\displaystyle
+\frac{\mu_1}{\mu_2}2(\mu_2-\mu_1)\left\vert Sym\nabla\mbox{\boldmath $u$}
_1-\frac{\nabla\cdot\mbox{\boldmath $u$}
_1}
{3}I_3\right\vert^2 dx\\
\\
\displaystyle
\le <(\Lambda_{\lambda_2,\,\mu_2}-\Lambda_{\lambda_1,\,\mu_1})\mbox{\boldmath $f$}
, \mbox{\boldmath $f$}>,
\end{array}
\tag {5.3}
$$
$$\begin{array}{c}
\displaystyle
<(\Lambda_{\lambda_2,\,\mu_2}-\Lambda_{\lambda_1,\,\mu_1})\mbox{\boldmath $f$}
, \mbox{\boldmath $f$}>\\
\\
\displaystyle
\le\int_{\Omega}\frac{3(\lambda_2-\lambda_1)+2(\mu_2-\mu_1)}{3}
\vert\nabla\cdot\mbox{\boldmath $u$}
_1\vert^2\\
\\
\displaystyle
+2(\mu_2-\mu_1)\left\vert Sym\nabla\mbox{\boldmath $u$}
_1
-\frac{\nabla\cdot\mbox{\boldmath $u$}
_1}{3}I_3\right\vert^2 dx.
\end{array}
\tag {5.4}
$$
\endproclaim

{\it\noindent Remark.}
These inequalities are a simple consequence of the system of integral inequalities in \cite{IK1},
the identities (3.8) and (3.9) in \cite{IK1} and
the factorization of symmetric matrix $B$:
$$
B=\frac{\text{Trace} B}{3}I_3 + \left(B-\frac{\text{Trace} B}{3}I_3\right).
$$
However, for reader's convenience, we present a direct proof.

{\it\noindent Proof.}
It is easy to see that
$$\begin{array}{c}
\displaystyle
<(\Lambda_{\lambda_2,\,\mu_2}-\Lambda_{\lambda_1,\,\mu_1})\mbox{\boldmath $f$}
, \mbox{\boldmath $f$}>
=\int_{\Omega}\lambda_1\vert\nabla\cdot(\mbox{\boldmath $u$}
_1-\mbox{\boldmath $u$}
_2)\vert^2
+2\mu_1\vert Sym\nabla(\mbox{\boldmath $u$}
_1-\mbox{\boldmath $u$}
_2)\vert^2\\
\\
\displaystyle
+(\lambda_2-\lambda_1)\vert\nabla\cdot\mbox{\boldmath $u$}
_2\vert^2
+2(\mu_2-\mu_1)\vert Sym\nabla\mbox{\boldmath $u$}
_2\vert^2 dx.
\end{array}
\tag {5.5}
$$
Set
$$\displaystyle
B_j=Sym\nabla\mbox{\boldmath $u$}
_j-\frac{\nabla\cdot\mbox{\boldmath $u$}
_j}{3}I_3.
$$
Note that $B_j\cdot I_3=\text{Trace}\, B_j=0$.
Since for any $\alpha, \beta$ and $3\times 3$ matrix $A$
$$\begin{array}{c}
\displaystyle
\alpha\vert\text{Trace}\,A\vert^2+2\beta\vert Sym A\vert^2\\
\\
\displaystyle
=\frac{3\alpha+2\beta}{3}\vert\text{Trace}\,A\vert^2
+2\beta\left\vert Sym A-\frac{\text{Trace}\, A}{3}I_3\right\vert^2,
\end{array}
$$
(5.5) can be rewritten as
$$\begin{array}{c}
\displaystyle
<(\Lambda_{\lambda_2,\,\mu_2}-\Lambda_{\lambda_1,\mu_1})\mbox{\boldmath $f$}
, \mbox{\boldmath $f$}
>\\
\\
\displaystyle
=
\int_{\Omega}\frac{3\lambda_1+2\mu_1}{3}\vert\nabla\cdot(\mbox{\boldmath $u$}
_1-\mbox{\boldmath $u$}
_2)\vert^2
+2\mu_1\vert B_1-B_2\vert^2\\
\\
\displaystyle
+\frac{3(\lambda_2-\lambda_1)+2(\mu_2-\mu_1)}{3}\vert\nabla\cdot\mbox{\boldmath $u$}
_2\vert^2
+2(\mu_2-\mu_1)\vert B_2\vert^2.
\end{array}
\tag {5.6}
$$
A combination of (E) and (5.6) gives
$$\begin{array}{c}
\displaystyle
<(\Lambda_{\lambda_2,\,\mu_2}-\Lambda_{\lambda_1,\,\mu_1})\mbox{\boldmath $f$}
, \mbox{\boldmath $f$}
>\\
\\
\displaystyle
\ge
\int_{\Omega}
\frac{3(\lambda_2-\lambda_1)+2(\mu_2-\mu_1)}{3}\vert\nabla\cdot\mbox{\boldmath $u$}
_2\vert^2
+2(\mu_2-\mu_1)\vert B_2\vert^2 dx.
\end{array}
$$
Replacing $2$ with $1$ and $1$ with $2$, we get (5.4).
On the other hand, from (E) we get
$$\begin{array}{c}
\displaystyle
\frac{3\lambda_1+2\mu_1}{3}\vert\nabla\cdot(\mbox{\boldmath $u$}
_1-\mbox{\boldmath $u$}
_2)\vert^2
+2\mu_1\vert B_1-B_2\vert^2\\
\\
\displaystyle
+\frac{3(\lambda_2-\lambda_1)+2(\mu_2-\mu_1)}{3}\vert\nabla\cdot\mbox{\boldmath $u$}
_2\vert^2
+2(\mu_2-\mu_1)\vert B_2\vert^2\\
\\
\displaystyle
=\frac{3\lambda_2+2\mu_2}{3}\vert\nabla\cdot\mbox{\boldmath $u$}
_2\vert^2
-\frac{2(3\lambda_1+2\mu_1)}{3}\nabla\cdot\mbox{\boldmath $u$}
_1\nabla\cdot\mbox{\boldmath $u$}
_2\\
\\
\displaystyle
+2\mu_2\vert B_2\vert^2-4\mu_1 B_1\cdot B_2\\
\\
\displaystyle
+\frac{3\lambda_1+2\mu_1}{3}\vert\nabla\cdot\mbox{\boldmath $u$}
_1\vert^2
+2\mu_1\vert B_1\vert^2\\
\\
\displaystyle
=\left\vert\sqrt{\frac{3\lambda_2+2\mu_2}{3}}\nabla\cdot\mbox{\boldmath $u$}
_2
-\sqrt{\frac{3}{3\lambda_2+2\mu_2}}\frac{3\lambda_1+2\mu_1}{3}\nabla\cdot\mbox{\boldmath $u$}
_1\right\vert^2\\
\\
\displaystyle
+\left\vert\sqrt{2\mu_2}B_2-\frac{2\mu_1}{\sqrt{2\mu_2}} B_1\right\vert^2\\
\\
\displaystyle
+\left\{\frac{3\lambda_1+2\mu_1}{3}-\frac{3}{3\lambda_2+2\mu_2}\left(\frac{3\lambda_1+2\mu_1}{3}\right)^2\right\}
\vert\nabla\cdot\mbox{\boldmath $u$}
_1\vert^2\\
\\
\displaystyle
+\left\{2\mu_1-\frac{(2\mu_1)^2}{2\mu_2}\right\}\vert B_1\vert^2\\
\\
\displaystyle
\ge
\frac{3\lambda_1+2\mu_1}{3(3\lambda_2+2\mu_2)}
\{3(\lambda_2-\lambda_1)+2(\mu_2-\mu_1)\}\vert\nabla\cdot\mbox{\boldmath $u$}
_1\vert^2
+\frac{2\mu_1(\mu_2-\mu_1)}{\mu_2}\vert B_1\vert^2
\end{array}
$$
Combining this with (5.6), we immediately get (5.3).  Similarly, we obtain (5.1) and (5.2), too.

\noindent
$\Box$

We construct for each $x\in\Omega$ two kinds of singular solutions of the equation
${\cal L}_{\lambda_0,\,\mu_0}\mbox{\boldmath $u$}
=0$ in $\Omega\setminus\{x\}$.

Let $G(\cdot)$ denote the standard fundamental solution for $-\triangle$:
$$\displaystyle
G(z)=\frac{1}{4\pi\vert z\vert}.
$$
The lemma below is a result of speciality of "isotropic" and the proof is given in Section 6.

\proclaim{\noindent
Proposition 5.2.}
Assume that
$$
\text{$\partial\Omega$ is $C^{2,1}$}
\tag {5.7}
$$
and that
$$
\text{$\lambda_0\in C^{0,1}(\overline\Omega)$, $\mu_0\in C^{2,1}(\overline\Omega)$.}
\tag {5.8}
$$
There exists a family $(\mbox{\boldmath $u$}
^0_x\in\{H^1_{\text{loc}}(\Omega\setminus\{x\})\}^3)_{x\in\Omega}$
such that
$$\displaystyle
\left(\mbox{\boldmath $u$}
^0_x-\nabla G(\cdot- x)+\frac{G(\cdot- x)}{\lambda_0(x)+2\mu_0(x)}
\left\{I_3-\frac{\cdot-x}{\vert \cdot-x\vert}\otimes\frac{\cdot-x}{\vert\cdot-x\vert}\right\}
\nabla\mu_0(x)\right)_{x\in\Omega}
$$
is bounded in $\{H^1(\Omega)\}^3$ and
$$
\text{${\cal L}_{\lambda_0,\,\mu_0}\mbox{\boldmath $u$}
^0_x=0$ in $\Omega\setminus\{x\}$.}
$$
\endproclaim

{\it\noindent Remark.}
$\nabla G(\cdot-x)$ satisfies the equation
$$\displaystyle
{\cal L}_{\lambda_0(x),\,\mu_0(x)}\mbox{\boldmath $v$}
+(\lambda_0(x)+2\mu_0(x))
\nabla\delta(\cdot-x)=0\,\,\text{in}\, \Bbb R^3
$$
and its divergence vanishes in $\Bbb R^3\setminus\{x\}$.

Let $\mbox{\boldmath $E$}
_0(\cdot;x)$ denote the standard fundamental solution(Kelvin matrix) for the operator
${\cal L}_{\lambda_0(x),\,\mu_0(x)}$ (see for example \cite{K}):
$$\displaystyle
\mbox{\boldmath $E$}
_0(z;x)=\frac{1}{8\pi}\left(\frac{1}{\mu_0(x)}
+\frac{1}{\lambda_0(x)+2\mu_0(x)}\right)\frac{I_3}{\vert z\vert}
+\frac{1}{8\pi}\left(\frac{1}{\mu_0(x)}-\frac{1}{\lambda_0(x)+2\mu_0(x)}\right)
\frac{z\otimes z}{\vert z\vert^3}.
$$

\proclaim{\noindent Proposition 5.3.}
Assume that $\lambda_0, \mu_0\in C^{0,1}(\overline\Omega)$.
There exists a family $(\mbox{\boldmath $E$}
^0_x(\cdot))_{x\in\Omega}$ in
$\cap_{1\le p<3}L^p(\Omega, M_3(\Bbb R))$ such that
$(\mbox{\boldmath $E$}
^0_x(\cdot)-\mbox{\boldmath $E$}
_0(\cdot-x;x))_{x\in\Omega}$
is bounded in $H^1(\Omega, M_3(\Bbb R))$
and for each constant vector $\mbox{\boldmath $b$}$
$$\displaystyle
{\cal L}_{\lambda_0,\,\mu_0}(\mbox{\boldmath $E$}
^0_x(\cdot)\mbox{\boldmath $b$}
)+\delta(\cdot-x)\mbox{\boldmath $b$}
=0
\,\,\text{in}\, \Omega.
$$
\endproclaim

The proof of this proposition proceeds  along the same lines with that of Proposition 2.1 and we may omit it.
We note that $\mbox{\boldmath $E$}
_0(\cdot-x;x)$ and $G(\cdot-x)$ are not independent of each other.
We will use the formula
$$\displaystyle
\nabla\cdot\{\mbox{\boldmath $E$}
_0(\cdot-x;x)\mbox{\boldmath $b$}
\}
=\frac{1}{\lambda_0(x)+2\mu_0(x)}\nabla G(\cdot-x)\cdot\mbox{\boldmath $b$},
\tag {5.9}
$$
where $\mbox{\boldmath $b$}$ is a constant vector.  That can be easily checked by direct computation.

\proclaim{\noindent Proposition 5.4(Selection of boundary data).}
Assume that $\partial\Omega$ is $C^3$ and that
$$\displaystyle
\lambda_0\in C^2(\overline\Omega), \mu_0\in C^3(\overline\Omega).
$$
Let $c$ be a needle.
\par
(i)  For each $t\in]0, 1[$ there exists a sequence
$\{\mbox{\boldmath $f$}
_n(\cdot ; c(t))\}$ of ${\cal D}(\Gamma)^3$  such that the weak solution
$\mbox{\boldmath $u$}_n$  of
$$\begin{array}{c}
\displaystyle
{\cal L}_{\lambda_0,\,\mu_0}\mbox{\boldmath $u$}
=0\quad\text{in}\, \Omega,\\
\\
\displaystyle
\mbox{\boldmath $u$}
\vert_{\partial\Omega}=\mbox{\boldmath $f$}
_n(\cdot ;c(t))
\end{array}
$$
converges to $\mbox{\boldmath $u$}
_{c(t)}^0(\cdot)$ in $\{H^1_{\text{loc}}(\Omega\setminus\{c(t')\vert
0<t'\le t \})\}^3$ as $n\longrightarrow\infty$;
\par
(ii)  Let $\{\mbox{\boldmath $e$}
_1, \mbox{\boldmath $e$}
_2, \mbox{\boldmath $e$}
_3\}$ be the standard orthonormal basis for $\Bbb R^3$.
For each $t\in]0, 1[$ and $j=1, 2, 3$ there exists a sequence
$\{\mbox{\boldmath $g$}
_{n,j}(\cdot ; c(t))\}$ of ${\cal D}(\Gamma)^3$ such that
the weak solution  $\mbox{\boldmath $v$}_n$  of
$$\begin{array}{c}
\displaystyle
{\cal L}_{\lambda_0,\,\mu_0}\mbox{\boldmath $v$}
=0\quad\text{in}\, \Omega,\\
\\
\displaystyle
\mbox{\boldmath $v$}
\vert_{\partial\Omega}=\mbox{\boldmath $g$}
_{n,j}(\cdot ;c(t))
\end{array}
$$
converges to $\mbox{\boldmath $E$}
^0_{c(t)}(\,\cdot\,)\mbox{\boldmath $e$}
_j$ in $\{H^1_{\text{loc}}(\Omega\setminus\{c(t')\vert
0<t'\le t \})\}^3$ as $n\longrightarrow\infty$.
\endproclaim

{\it\noindent Proof.}
Apply the Runge approximation property for the equation
${\cal L}_{\lambda_0, \mu_0}\mbox{\boldmath $u$}
=0$.  It is easily proved by using
the unique continuation property established in \cite{AITY}.

\noindent
$\Box$

{\bf\noindent Definition($T^f(c)$, $T^g(c)$, $I^f(t,c)$ and $I^g(t,c)$).}

(i)  Denote by $T^f(c)$ the set of all $s\in]0, 1[$ such that
$$\displaystyle
I^f(t,c)\equiv\lim_{n\rightarrow\infty}
 <(\Lambda_{\lambda,\,\mu}-\Lambda_{\lambda_0,\,\mu_0})\mbox{\boldmath $f$}
_n(\cdot;c(t))\,,\mbox{\boldmath $f$}
_n(\cdot;c(t))>
$$
exists for $0<t<s$ and $\sup_{0<t<s}\vert I^f(t,c)\vert<\infty$.

(ii)
Denote by $T^g(c)$ the set of all $s\in]0, 1[$ such that
$$\displaystyle
I^g(t,c)\equiv\lim_{n\rightarrow\infty}
\sum_{1\le j\le 3}
<(\Lambda_{\lambda, \mu}-\Lambda_{\lambda_0, \mu_0})
\mbox{\boldmath $g$}_{n,j}(\cdot;c(t))\,,\mbox{\boldmath $g$}
_{n,j}(\cdot;c(t))>
$$
exists for $0<t<s$ and $\sup_{0<t<s}\vert I^g(t,c)\vert<\infty$.

A combination of Propositions 5.1 and 5.4 yields

\proclaim{\noindent
Proposition 5.5.}
Let $c$ be a needle.  If $0<t<t(c;D)$, $I^f(t,c)$ and $I^g(t,c)$
exist and it holds that
$$\begin{array}{c}
\displaystyle
\int_D\frac{3\lambda_0+2\mu_0}{3(3\lambda+2\mu)}
\{3(\lambda-\lambda_0)+2(\mu-\mu_0)\}\vert\nabla\cdot\mbox{\boldmath $u$}
^0_{c(t)}\vert^2\\
\\
\displaystyle
+\frac{\mu_0}{\mu}2(\mu-\mu_0)\left\vert Sym\nabla\mbox{\boldmath $u$}
^0_{c(t)}-\frac{\nabla\cdot
\mbox{\boldmath $u$}
^0_{c(t)}}{3}I_3\right\vert^2 dy\le  I^f(t,c),
\end{array}
\tag {5.10}
$$
$$\begin{array}{c}
\displaystyle
I^f(t,c)\le\int_D\frac{3(\lambda-\lambda_0)+2(\mu-\mu_0)}{3}
\vert\nabla\cdot\mbox{\boldmath $u$}
^0_{c(t)}\vert^2\\
\\
\displaystyle
+2(\mu-\mu_0)\left\vert Sym\nabla\mbox{\boldmath $u$}
^0_{c(t)}
-\frac{\nabla\cdot\mbox{\boldmath $u$}
^0_{c(t)}}{3}I_3\right\vert^2 dy,
\end{array}
\tag {5.11}
$$
and
$$
\begin{array}{c}
\displaystyle
\sum_{1\le j\le 3}\int_D\frac{3\lambda_0+2\mu_0}{3(3\lambda+2\mu)}
\{3(\lambda-\lambda_0)+2(\mu-\mu_0)\}\vert\nabla\cdot(\mbox{\boldmath $E$}
^0_{c(t)}\mbox{\boldmath $e$}
_j)
\vert^2\\
\\
\displaystyle
+\frac{\mu_0}{\mu}2(\mu-\mu_0)\left\vert Sym\nabla(\mbox{\boldmath $E$}
^0_{c(t)}\mbox{\boldmath $e$}
_j)-
\frac{\nabla\cdot(\mbox{\boldmath $E$}
^0_{c(t)}\mbox{\boldmath $e$}
_j)}{3}I_3\right\vert^2 dy\le  I^g(t,c),
\end{array}
\tag {5.12}
$$
$$\begin{array}{c}
\displaystyle
I^g(t,c)\le\sum_{1\le j\le 3}
\int_D\frac{3(\lambda-\lambda_0)+2(\mu-\mu_0)}{3}
\vert\nabla\cdot(\mbox{\boldmath $E$}
^0_{c(t)}\mbox{\boldmath $e$}
_j)\vert^2\\
\\
\displaystyle
+2(\mu-\mu_0)\left\vert Sym\nabla(\mbox{\boldmath $E$}
^0_{c(t)}\mbox{\boldmath $e$}
_j)
-\frac{\nabla\cdot(\mbox{\boldmath $E$}
^0_{c(t)}\mbox{\boldmath $e$}
_j)}{3}I_3\right\vert^2 dy.
\end{array}
\tag {5.13}
$$
\endproclaim

The proof is same as that of Proposition 2.4 and we may omit it.

Theorem D is a consequence of

\proclaim{\noindent Proposition 5.6(Reconstruction of impact parameter).}
\newline{For any needle $c$} it holds that
$$
T^f(c)\cap T^g(c)=]0, t(c;D)[
$$
and thus the formula
$$
t(c;D)=\sup T^f(c)\cap T^g(c),
$$
is valid.  $\partial D$ has the form
$$
\partial D=\{c(\sup T^f(c)\cap T^g(c))\,\vert\, \text{$c$ is a needle and satisfies
$\sup T^f(c)\cap T^g(c)<1$}\}.
$$
\endproclaim

{\it\noindent Remark.}  ${\cal D}(\lambda_0, \mu_0, \Gamma)$ in Theorem D is
$$
\{\mbox{\boldmath $f$}
_n(\cdot;c(t))\,, \mbox{\boldmath $g$}
_{n,j}(\cdot;c(t))\,,
j=1,2,3,\,n=1,\cdots\,\vert\,\text{$c$ is a needle and $0<t<1$}\}.
$$

{\it\noindent Proof.}
First from Propositions 5.2, 5.3 and 5.5 we get $]0, t(c;D)[\subset T^f(c)\cap T^g(c)$.
Since $T^f(c)\cap T^g(c)\subset ]0, 1[$, if $t(c;D)=1$, we know that $T^f(c)\cap T^g(c)=]0, 1[=]0, t(c;D)[$.

So the problem is the case where $t(c;D)<1$.  Set $a=c(t(c;D))(\in\partial D)$.
We assume that $]0, t(c;D)[$ does not coincide with $T^f(c)\cap T^g(c)$.  Then there exists
$s\in T^f(c)\cap T^g(c)$
such that $t(c;D)\le s$.  Then from the definitions of $T^f(c)$, $T^g(c)$ we get
$$
\text{$\exists C>0\,\forall t\in]0, t(c;D)[\,\,\,
\vert I^f(t,c)\vert<C$ and $\vert I^g(t,c)\vert<C$.}
\tag {5.14}
$$
{\bf  Case 1. $\mu_D(a)\not=\mu_0(a)$}
\newline{We consider first} the case when $\mu_D(a)>\mu_0(a)$.
In what follows, $C_1, C_2, \cdots$ denote positive constants independent of $x\in\Omega$.
Then we may assume that
$$
\exists r>0\,\exists C_1>0\,\forall y\in B(a, 2r)\cap D\,\,\mu(y)-\mu_0(y)\ge C_1.
\tag {5.15}
$$
take $t_1\in ]0, t(c;D)[$ such that $c(t)\in B(a, r)$ for $t\in ]t_1, t(c;D)[$.
Then  we know that $\vert y-c(t)\vert\ge r$ for $y\in\Omega\setminus B(a, 2r)$.
From Proposition 5.2 we get
$$\displaystyle
\exists C_2>0\,\int_{D\setminus B(a,2r)}
\vert Sym\nabla\mbox{\boldmath $u$}
_{c(t)}^0\vert^2 dy <C_2.
\tag {5.16}
$$
On the other hand, since $c(t)\in \Omega\setminus\overline D$ for $t_1<t<t(c;D)$,
we deduce that
$$
\text{$\nabla\cdot\mbox{\boldmath $u$}
_{c(t)}^0=\nabla\cdot\{\mbox{\boldmath $u$}
_{c(t)}^0-
\nabla G(\cdot-c(t))\}$ in $D\cap B(a, 2r)$.}
$$
Thus we get
$$\displaystyle
\int_{D\cap B(a,2r)}\vert\nabla\cdot\mbox{\boldmath $u$}
_{c(t)}^0\vert^2 dy
\le C_3\int_{D\cap B(a,2r)}\frac{1}{\vert y-c(t)\vert^4} dy + C_4.
\tag {5.17}
$$
Noting that $\vert\nabla\cdot\mbox{\boldmath $u$}
_{c(t)}^0\vert\le\sqrt{3}\vert Sym\nabla\mbox{\boldmath $u$}
_{c(t)}^0\vert$,
from (5.10), (5.16) and (5.17) we get
$$\begin{array}{c}
\displaystyle
\int_{D\cap B(a, 4r)}\frac{\mu}{\mu_0}2(\mu-\mu_0)
\left\vert Sym\nabla\mbox{\boldmath $u$}
_{c(t)}^0-\frac{\nabla\cdot\mbox{\boldmath $u$}
_{c(t)}^0}{3}I_3\right\vert^2dy\\
\\
\displaystyle
\le I^f(t,c)+C_5\int_{D\cap B(a,2r)}\frac{1}{\vert y-c(t)\vert^4} dy + C_6.
\end{array}
\tag {5.18}
$$
Since
$$\displaystyle
\left\vert Sym\nabla\mbox{\boldmath $u$}
_{c(t)}^0-\frac{\nabla\cdot\mbox{\boldmath $u$}
_{c(t)}}{3}I_3\right\vert^2
\ge \frac{1}{2}\vert Sym\nabla\mbox{\boldmath $u$}
_{c(t)}^0\vert^2
-\left\vert\frac{\nabla\cdot\mbox{\boldmath $u$}
_{c(t)}^0}{3}I_3\right\vert^2,
$$
it follows from (5.15), (5.17) and (5.18) that
$$\begin{array}{c}
\displaystyle
\int_{D\cap B(a,2r)}\vert Sym\nabla\mbox{\boldmath $u$}
_{c(t)}^0\vert^2dy\\
\\
\displaystyle
\le C_7 I^f(t,c)+C_8\int_{D\cap B(a,2r)}\frac{1}{\vert y-c(t)\vert^4} dy + C_9.
\end{array}
\tag {5.19}
$$
A combination of the triangle inequality and Proposition 5.2 yields
$$\begin{array}{c}
\displaystyle
\int_{D\cap B(a,2r)}\vert Sym\nabla\mbox{\boldmath $u$}
_{c(t)}^0\vert^2dy\\
\\
\displaystyle
\ge C_{10}\int_{D\cap B(a,2r)}\frac{1}{\vert y-c(t)\vert^6}dy
- C_{11}\int_{D\cap B(a,2r)}\frac{1}{\vert y-c(t)\vert^4}dy-C_{12}.
\end{array}
\tag {5.20}
$$
From (5.19),  (5.20) we get
$$\begin{array}{c}
\displaystyle
\int_{D\cap B(a,2r)}\frac{1}{\vert y-c(t)\vert^6}dy\\
\\
\displaystyle
\le C_7 I^f(t,c)+C_{13}\int_{D\cap B(a,2r)}\frac{1}{\vert y-c(t)\vert^4} dy + C_{14}.
\end{array}
\tag {5.21}
$$
The Young inequality yields
$$\begin{array}{c}
\displaystyle
\int_{D\cap B(a,2r)}\frac{1}{\vert y-c(t)\vert^4}dy\\
\\
\displaystyle
\le \frac{1}{3\epsilon ^3}\vert D\cap B(a,2r)\vert+\frac{2\epsilon ^{3/2}}{3}\int_{D\cap B(a,2r)}
\frac{1}{\vert y-c(t)\vert^6}dy,
\,\forall\epsilon>0.
\end{array}
\tag {5.22}
$$
If we take a sufficiently small $\epsilon$, from (5.21) and (5.22) we get
$$\displaystyle
\displaystyle
C_{15}\int_{D\cap B(a,2r)}\frac{1}{\vert y-c(t)\vert^6}dy-C_{16}
\le I^f(t,c).
\tag {5.23}
$$
Since $\partial D$ is Lipschitz, it is easy to see that
$$\displaystyle
\lim_{t\uparrow t(c;D)}\int_{D\cap B(a,2r)}\frac{1}{\vert y-c(t)\vert^6}dy=\infty
$$
and (5.23) therefore yields $\lim_{t\uparrow t(c;D)}I^f(t,c)=\infty$.
This  is a contradiction to (5.14).

We now consider the case when $\mu_D(a)<\mu_0(a)$.
Using (5.11), we get
\newline{$\lim_{t\uparrow t(c;D)}I^f(t,c)=-\infty$}
and a contradiction to (5.14).  The proof proceeds along the same lines.
\newline{{\bf Case 2. $\mu_D(a)=\mu_0(a)$}}
\newline{From (v) of (EI)}
we know that $\lambda_D(a)\not=\lambda_0(a)$.  Let $\epsilon>0$.
Then there exists a positive constant $C_a$ such that if we take a
small $r>0$ , we may assume that
$$\displaystyle
\vert\mu(y)-\mu_0(y)\vert<\epsilon, \forall y\in D\cap B(a, 2r)
\tag {5.24}
$$
and
$$\displaystyle
\lambda(y)-\lambda_0(y)\ge C_a, \forall y\in D\cap B(a, 2r)
\tag {5.25}
$$
or
$$\displaystyle
\lambda(y)-\lambda_0(y)\le - C_a, \forall y\in D\cap B(a, 2r).
\tag {5.26}
$$
Let us consider the case when (5.24) and (5.25) are satisfied.
Take $t_1\in]0, t(c;D)[$ such that $c(t)\in B(a,r)$ for $t\in]t_1, t(c;D)[$.
Then  we know that $\vert y-c(t)\vert\ge r$ for $y\in\Omega\setminus B(a,2r)$.
From Proposition 5.3 we get
$$\begin{array}{c}
\displaystyle
\vert\sum_{1\le j\le 3}\int_{D\setminus B(a,2r)}\frac{3\lambda_0+2\mu_0}{3(3\lambda+2\mu)}
\{3(\lambda-\lambda_0)+2(\mu-\mu_0)\}\vert\nabla\cdot(\mbox{\boldmath $E$}
^0_{c(t)}\mbox{\boldmath $e$}
_j)
\vert^2\\
\\
\displaystyle
+\frac{\mu_0}{\mu}2(\mu-\mu_0)\left\vert Sym\nabla(\mbox{\boldmath $E$}
^0_{c(t)}\mbox{\boldmath $e$}
_j)-
\frac{\nabla\cdot(\mbox{\boldmath $E$}
^0_{c(t)}\mbox{\boldmath $e$}
_j)}{3}I_3\right\vert^2 dy\vert \le C(r),
\end{array}
$$
and therefore (5.12) yields
$$\begin{array}{c}
\displaystyle
\sum_{1\le j\le 3}\int_{D\cap B(a,2r)}\frac{3\lambda_0+2\mu_0}{3(3\lambda+2\mu)}
\{3(\lambda-\lambda_0)+2(\mu-\mu_0)\}\vert\nabla\cdot(\mbox{\boldmath $E$}
^0_{c(t)}\mbox{\boldmath $e$}
_j)
\vert^2\\
\\
\displaystyle
+\frac{\mu_0}{\mu}2(\mu-\mu_0)\left\vert Sym\nabla(\mbox{\boldmath $E$}
^0_{c(t)}\mbox{\boldmath $e$}
_j)-
\frac{\nabla\cdot(\mbox{\boldmath $E$}
^0_{c(t)}\mbox{\boldmath $e$}
_j)}{3}I_3\right\vert^2 dy-C(r)\\
\\
\displaystyle
\le I^g(t,c).
\end{array}
\tag {5.27}
$$
A combination of (5.24), Proposition 5.3 and the triangle inequality yields
$$\begin{array}{c}
\displaystyle
\vert\sum_{1\le j\le 3}\int_{D\cap B(a,2r)}\frac{3\lambda_0+2\mu_0}{3(3\lambda+2\mu)}
2(\mu-\mu_0)\vert\nabla\cdot(\mbox{\boldmath $E$}
^0_{c(t)}\mbox{\boldmath $e$}
_j)
\vert^2\\
\\
\displaystyle
+\frac{\mu_0}{\mu}2(\mu-\mu_0)\left\vert Sym\nabla(\mbox{\boldmath $E$}
^0_{c(t)}\mbox{\boldmath $e$}
_j)-
\frac{\nabla\cdot(\mbox{\boldmath $E$}
^0_{c(t)}\mbox{\boldmath $e$}
_j)}{3}I_3\right\vert^2 dy\vert\\
\\
\displaystyle
\le C_{17}\epsilon \left(1+\int_{D\cap B(a,2r)}\frac{1}{\vert y-c(t)\vert^4}dy\right).
\end{array}
\tag {5.28}
$$
From (5.9), (5.25), Proposition 5.3 and the triangle inequality we get
$$\begin{array}{c}
\displaystyle
\sum_{1\le j\le 3}\int_{D\cap B(a,2r)}\frac{3\lambda_0+2\mu_0}{3(3\lambda+2\mu)}
3(\lambda-\lambda_0)\vert\nabla\cdot(\mbox{\boldmath $E$}
^0_{c(t)}\mbox{\boldmath $e$}
_j)
\vert^2 dy\\
\\
\displaystyle
\ge C_{18}\int_{D\cap B(a,2r)}\vert\nabla G(y-c(t))\vert^2dy
-C_{19}.
\end{array}
\tag {5.29}
$$
A combination of (5.27) to (5.29) yields
$$\displaystyle
(C_{18}-C_{17}\epsilon)\int_{D\cap B(a,2r)}\frac{1}{\vert y-c(t)\vert^4}dy
-C_{17}\epsilon-C_{19}-C(r)\le I^g(t,c).
\tag {5.30}
$$
Notice that $C_{17}$ and $C_{18}$ are independent of $r$ and $\epsilon$.
So if we take $\epsilon$ in such a way that $C_{18}-C_{17}\epsilon>0$ in advance,
we get $\lim_{t\uparrow t(c;D)}I^g(t,c)=\infty$.  This is a contradiction to (5.14).
We can apply a similar argument to (5.13) in the case when (5.24) and (5.26) are satisfied
and get $\lim_{t\uparrow t(c;D)}I^g(t,c)=-\infty$.
This completes the proof of Proposition 5.6.

\noindent
$\Box$

\section{Construction of singular solutions}

{\it\noindent Proof of Proposition 2.1.}
Set
$$
s=
\left\{
\begin{array}{lr}
\displaystyle
1 & \text{if $n=3$}\\
\\
\displaystyle
\theta & \text{ if $n=2$}.
\end{array}
\right.
$$
We construct $G_x^0(\cdot)$ in the form
$$
G_x^0(\cdot)=G_0(\cdot-x;x)+\epsilon(\cdot;x).
$$
Define a functional $f_x$ by the formula
$$
f_x(\varphi)=\int_{\Omega}
(\gamma_0(x)-\gamma_0(y))
\nabla G_0(y-x;x)\cdot\nabla\varphi(y)dy,
\varphi\in C_0^{\infty}(\Omega).
$$
By the assumption, we have
$$
\vert\gamma_0(x)-\gamma_0(y)\vert\le M\vert x-y\vert^{s}.
$$
Since
$$
\vert\nabla G_0(y-x;x)\vert\le\frac{C}{\vert y-x\vert^{n-1}},
$$
we get
$$\begin{array}{c}
\displaystyle
\vert f_x(\varphi)\vert\le
\int_{\Omega}
\vert\gamma_0(x)-\gamma_0(y)\vert
\vert\nabla G_0(y-x;x)\vert\vert\nabla\varphi(y)\vert dy\\
\\
\displaystyle
\le C\left\{\int_{\Omega}
\frac{dy}{\vert y-x\vert^{2(n-1-s)}}\right\}^{1/2}\vert\varphi\vert_{1,\Omega}.
\end{array}
$$
Therefore $f_x$ is in $(H^1_0(\Omega))^*$ since $2(n-1-s)<n$.
Let $\epsilon(\cdot;x)\in H^1_0(\Omega)$ be the weak solution of
$$\begin{array}{c}
\displaystyle
\nabla\cdot\gamma_0\nabla\epsilon(\cdot;x)=-f_x\,\,\text{in}\,\Omega,\\
\\
\displaystyle
\epsilon(\cdot;x)\vert_{\partial\Omega}=0.
\end{array}
$$
$\epsilon(\cdot;x)$ satisfies
$$
\displaystyle
\vert\epsilon(\cdot;x)\vert_{1,\Omega}
\le C\left\{\int_{\Omega}\frac{dy}{\vert y-x\vert^{2(n-1-s)}}\right\}^{1/2}.
$$
For any $\varphi\in C_0^{\infty}(\Omega)$
$$\begin{array}{c}
\displaystyle
\int_{\Omega}\gamma_0(y)\nabla G_x^0(y)\cdot\nabla\varphi(y)dy\\
\\
\displaystyle
=\int_{\Omega}\gamma_0(y)\nabla G_0(y-x;x)\cdot\nabla\varphi(y)
+\gamma_0(y)\nabla\epsilon(y;x)\cdot\nabla\varphi(y)dy\\
\\
\displaystyle
=\int_{\Omega}\gamma_0(x)\nabla G_0(y-x;x)\cdot\nabla\varphi(y) dy\\
\\
\displaystyle
=\varphi(x).
\end{array}
$$
$G_x^0(\cdot)$ satisfies
$$\displaystyle
\vert G_x^0(\cdot)-G_0(\cdot-x;x)\vert_{1,\Omega}
\le C\left\{\int_{\Omega}\frac{dy}{\vert y-x\vert^{2(n-1-s)}}\right\}^{1/2}.
$$
From this inequality, we know that
$(G_x^0(\cdot)-G_0(\cdot-x;x))_{x\in\Omega}$ is bounded in $H^1(\Omega)$.

\noindent
$\Box$

{\it\noindent Proof of Proposition 5.2.}

We construct $\mbox{\boldmath $u$}
_{c(t)}^0$ in the form
$$\displaystyle
\mbox{\boldmath $u$}
_{c(t)}^0=\nabla G(\cdot-x)+{\cal E}(\cdot;x).
$$
Since
$$
\text{${\cal L}_{\lambda_0(x),\,\mu_0(x)}\nabla G(\cdot-x)=0$ in $\Omega\setminus\{x\}$,}
$$
it suffices to construct ${\cal E}(\cdot;x)$ such that
$$
\displaystyle
{\cal L}_{\lambda_0,\,\mu_0}{\cal E}(\cdot;x)=\{{\cal L}_{\lambda_0(x),\,\mu_0(x)}
-{\cal L}_{\lambda_0,\,\mu_0}\}\nabla G(\cdot-x)\,\,\text{in}\, \Omega\setminus\{x\}.
\tag {6.1}
$$
In fact, this is divided into two steps.

First we construct ${\cal E}_1(\cdot;x)$ such that
$$
\displaystyle
{\cal L}_{\lambda_0(x),\,\mu_0(x)}
{\cal E}_1(\cdot;x)=\{{\cal L}_{\lambda_0(x),\,\mu_0(x)}
- {\cal L}_{\lambda_0,\,\mu_0}\}\nabla G(\cdot-x)\,\,\text{in}\, \Omega\setminus\{x\}.
\tag {6.2}
$$

Second we construct ${\cal E}_2(\cdot;x)$
such that
$$
\displaystyle
{\cal L}_{\lambda_0,\,\mu_0}{\cal E}_2(\cdot;x)
=\{{\cal L}_{\lambda_0(x),\,\mu_0(x)}
-{\cal L}_{\lambda_0,\,\mu_0}\}{\cal E}_1(\cdot;x)\,\,\text{in}\,\Omega.
\tag {6.3}
$$
Then
$$\displaystyle
{\cal E}={\cal E}_1+{\cal E}_2
$$
satisfies (6.1).
Let us start with the explanation of

{\bf\noindent Construction of ${\cal E}_1(\cdot;x)$}

This is based on

\proclaim {\noindent Claim 1.}
Assume that $w$ and ${\cal E}_1^0$ satisfy
$$\displaystyle
\triangle w=-\nabla G(\cdot-x)\cdot\nabla\mu_0(\cdot)\,\,\text{in}\, \Omega\setminus\{x\}
\tag {6.4}
$$
and
$$\displaystyle
{\cal L}_{\lambda_0(x),\,\mu_0(x)}{\cal E}_1^0=2\{\nabla\nabla\mu_0(\,\cdot\,)\}
\nabla G(\cdot-x)\,\,\text{in}\, \Omega\setminus\{x\}.
\tag {6.5}
$$
Then
$$\displaystyle
{\cal E}_1=\frac{2\nabla w}{\lambda_0(x)+2\mu_0(x)}
+{\cal E}_1^0
\tag {6.6}
$$
satisfies (6.2).
\endproclaim

This claim is easily checked if one knows
$$
{\cal L}_{\lambda_0(x), \mu_0(x)}(\nabla w)
= \{\lambda_0(x)+2\mu_0(x)\}\nabla\triangle w
$$
and
$$\begin{array}{c}
\displaystyle
\{{\cal L}_{\lambda_0(x),\,\mu_0(x)}-{\cal L}_{\lambda_0,\,\mu_0}\}
\nabla G(\cdot-x)\\
\\
\displaystyle
=-2\nabla\{\nabla G(\cdot-x)\cdot\nabla\mu_0(\cdot)\}
+2\{\nabla\nabla\mu_0(\cdot)\}\nabla G(\cdot-x).
\end{array}
$$
Therefore we have to explain how to construct such $w$ and ${\cal E}^0_1$.

{\it Construction of $w$.}

We construct $w$ in the form
$$\displaystyle
w=-\nabla\cdot\xi^0-\nabla\cdot\xi^1+\eta
\tag {6.7}
$$
where
$$\displaystyle
\xi^0=\xi^0(\cdot;x)=\frac{\vert\cdot-x\vert}{8\pi}\nabla\mu_0(x);
\tag {6.8}
$$
$\xi^1=\xi^1(\cdot;x)\in \{H^1_0(\Omega)\}^3$ is the weak solution of
$$\displaystyle
\triangle\xi^1=G(\cdot-x)\{\nabla\mu_0(\cdot)-\nabla\mu_0(x))\,\,
\text{in}\, \Omega;
\tag {6.9}
$$
$\eta=\eta(\cdot;x)\in H^1_0(\Omega)$ is the weak solution of
$$
\displaystyle
\triangle\eta=G(\cdot-x)\triangle\mu_0(\cdot)\,\,\text{in}\,\Omega.
\tag {6.10}
$$
Notice that
$\xi^0$ satisfies the equation
$$\displaystyle
\triangle\xi^0=G(\cdot-x)\nabla\mu_0(x)\,\,\text{in}\,\Omega.
\tag {6.11}
$$
A combination of (6.9), (6.10) and (6.11) implies that $w$ given by (6.7)
satisfies (6.4).  From (5.7), (5.8) and the regularity theory we know
$(\xi^1(\cdot;x))_{x\in\Omega}$ and $(\eta(\cdot;x))_{x\in\Omega}$
are bounded in $\{H^3(\Omega)\}^3$ and  $H^2(\Omega)$, respectively.
Therefore we can conclude
$$
\text{$\displaystyle\left(\nabla w+\frac{1}{8\pi}\{\nabla\nabla\{\vert\cdot-x\vert\}\}\nabla\mu_0(x)\right)_{x\in\Omega}$
is bounded in $\{H^1(\Omega)\}^3$.}
\tag {6.12}
$$

{\it Construction of ${\cal E}^0_1$.}
Since
$$\displaystyle
\{\nabla\nabla\mu_0(\cdot)-\nabla\nabla\mu_0(x)\}\nabla G(\cdot-x)\in \{L^2(\Omega)\}^3,
$$
We can find the weak solution ${\cal E}_1^{0,1}={\cal E}_1^{0,1}(\cdot;x)
\in \{H^1_0(\Omega)\}^3$ of
$$
\displaystyle
{\cal L}_{\lambda_0(x),\,\mu_0(x)}{\cal E}_1^{0,1}=2
\{\nabla\nabla\mu_0(\cdot)-\nabla\nabla\mu_0(x)\}\nabla G(\cdot-x)\,\,\text{in}\, \Omega.
\tag {6.14}
$$
We construct ${\cal E}_1^0$ in the form
$$\displaystyle
{\cal E}_1^0={\cal E}_1^{0,0}+{\cal E}_1^{0,1},
\tag {6.15}
$$
where ${\cal E}_1^{0,0}={\cal E}_1^{0,0}(\cdot;x)$ solves
$$\displaystyle
{\cal L}_{\lambda_0(x),\,\mu_0(x)}{\cal E}_1^{0,0}=2
\{\nabla\nabla\mu_0(x)\}\nabla G(\cdot-x)\,\,\text{in}\, \Omega.
\tag {6.16}
$$
The construction of ${\cal E}_1^{0,0}$ is based on

\proclaim{\noindent Claim 2.}  Assume that $f$ and $\mbox{\boldmath $u$}$ satisfy
$$\displaystyle
(\lambda_0(x)+2\mu_0(x))\triangle f=-(\lambda_0(x)+\mu_0(x))\nabla\cdot\mbox{\boldmath $u$}.
\tag {6.17}
$$
Then the formula
$$\displaystyle
{\cal L}_{\lambda_0(x),\,\mu_0(x)}(\mbox{\boldmath $u$}+\nabla f)=\mu_0(x)\triangle
\mbox{\boldmath $u$},
$$
is valid.
\endproclaim

This is directly checked and we may omit the proof.
Thus first we solve
$$\displaystyle
\mu_0(x)\triangle\mbox{\boldmath $u$}=
2\{\nabla\nabla\mu_0(x)\}\nabla G(\cdot-x)\,\,\text{in}\,\Omega.
$$
Fortunately, we can find the explicit solution given by
$$
\mbox{\boldmath $u$}=\mbox{\boldmath $u$}(\cdot;x)=\left\{\frac{\nabla\nabla\mu_0(x)}{4\pi\mu_0(x)}\right\}
\frac{(\cdot-x)}{\vert\cdot-x\vert}.
$$
$(\mbox{\boldmath $u$}(\cdot;x))_{x\in\Omega}$ is bounded in $\{H^1(\Omega)\}^3$
and thus we get the unique weak solution $f=f(\cdot;x)\in H^1_0(\Omega)$
of (6.17).  From the regularity theory and (5.7) we know
$(f(\cdot;x))_{x\in\Omega}$ is bounded in $H^2(\Omega)$.
Therefore ${\cal E}_1^{0,0}$ given by
$$\displaystyle
{\cal E}_1^{0,0}=\mbox{\boldmath $u$}+\nabla f
$$
satisfies (6.16) and hence ${\cal E}_1^0$ given by (6.15) is in
$\{H^1(\Omega)\}^3$. The boundedness of $({\cal
E}_1^0(\cdot;x))_{x\in\Omega}$ in $\{H^1(\Omega)\}^3$ is clear.
This completes the construction of ${\cal E}_1$.

{\bf\noindent Construction of ${\cal E}_2(\cdot;x)$.}

From the
construction of ${\cal E}_1$ and (5.8) we get the functional
$$\displaystyle
\{{\cal L}_{\lambda_0(x),\,\mu_0(x)}-{\cal L}_{\lambda_0,\,\mu_0}\}
{\cal E}_1(\cdot;x)
$$
belongs to the dual space of $\{H^1(\Omega)\}^3$. Thus one can
find the weak solution ${\cal E}_2={\cal E}_2(\cdot;x)\in
\{H_0^1(\Omega)\}^3$ of (6.3).  The boundedness of $({\cal
E}_2(\cdot;x))_{x\in\Omega}$ in $\{H^1(\Omega)\}^3$ is clear. This
completes the construction of ${\cal E}_2$.

\noindent
$\Box$

\centerline{{\bf Acknowledgement}}
The author thanks the referee for several suggestions for
the improvement of the manuscript.
This research was partially supported by Grant-in-Aid for Scientific
Research(No.  11640151), Ministry of Education, Science and Culture, Japan.

\section{Appendix}

\par For reader's convenience, we give a proof of

\proclaim{\noindent The Runge approximation property.}   Let $\gamma_0\in C^{0,1}(\overline\Omega)(n\ge 3)$,
$\in L^{\infty}(\Omega)(n=2)$.  Let $\Gamma$ be a given nonempty open subset of $\partial\Omega$.
Let $U$ be an open subset of $\Omega$.  Assume that
$\overline U\subset\Omega$ and that $\Omega\setminus\overline U$
is connected.
Then for any $u$ satisfying
$\nabla\cdot\gamma_0\nabla u=0$ in an open neighbourhood of $\overline U$
there exists a sequence $(u_k)$ of $H^1(\Omega)$ functions such that
$$\begin{array}{c}
\displaystyle
\nabla\cdot\gamma_0\nabla u_k=0\,\,\text{in}\,\Omega,\\
\\
\displaystyle
\text{supp}\,(u_k\vert_{\partial\Omega})\subset\Gamma,\\
\\
\displaystyle
u_k\longrightarrow u\,\,\text{in}\,H^1(U).
\end{array}
$$
\endproclaim

{\it\noindent Remark.}
Notice that if $\Gamma=\partial\Omega$, this
is a usual Runge approximation theorem.  A relationship between
the uniqueness for the Cauchy problem and the Runge approximation
property are described in \cite{LP} on pp.761-763.

$\quad$

{\it\noindent Proof.}
Extend $\gamma_0$ outside $\Omega$ as $1$ if $n=2$,
as a uniformly positive Lipschitz function if $n\ge 3$.
Take an open ball $B$ with small radius centered at a point in $\Gamma$
such that $B\cap\partial\Omega$ is contained in $\Gamma$ and represented
by a graph of a Lipschitz function on $\Bbb R^{n-1}$.
Then modifying it, we can construct a domain
$\Omega_0$ with Lipschitz boundary
in such a way that $\Omega\subset\Omega_0$,
$\Omega\not=\Omega_0$, and $\partial\Omega\setminus(B\cap\partial\Omega)$
is contained in $\partial\Omega_0$.
Give $F\in (H^1_0(\Omega_0))^*$ and find $u\in H^1_0(\Omega_0)$ such that
$$\begin{array}{c}
\displaystyle
\nabla\cdot\gamma_0\nabla u=-F\,\,\text{in}\,\Omega_0,\\
\\
\displaystyle
u\vert_{\partial\Omega_0}=0.
\end{array}
$$
Set
$$\displaystyle
\mbox{\boldmath $G$}F=u.
$$
We show that:

Let $X$ denote the set of all $u\vert_U$ satisfying $u\in H^1$ in
an open neighbourhood of $\overline U$ and therein
$\nabla\cdot\gamma_0\nabla u=0$; let $Y$ denote the set of all
$\mbox{\boldmath $G$}F\vert_{U}$ with
$\text{supp}\,F\subset\Omega_0\setminus\overline\Omega$. Then $Y$
is dense in $X$ with respect to $H^1$-topology.

By the Hahn-Banach theorem, this is equivalent to the following statement.
Let $f\in (H^1(U))^*$.  If
$$
f(\mbox{\boldmath $G$}F\vert_U)=0
$$
for all $F\in (H^1_0(\Omega_0))^*$ with
$\text{supp}\,F\subset\Omega_0\setminus\overline\Omega$ then $f$ vanishes on $X$, too.
\par
So we prove it.
Give $\varphi\in H^1_0(\Omega_0)$.  Set
$$
\tilde{f}(\varphi)=f(\varphi\vert_U).
$$
Of course $\tilde{f}$ is in $(H^1_0(\Omega_0))^*$.
Consider $\mbox{\boldmath $G$}\tilde{f}$.  Since $\text{supp}\,\tilde{f}\subset\overline U$,
we know that
$$\displaystyle
\nabla\cdot\gamma_0\nabla(\mbox{\boldmath $G$}\tilde{f})=0\,\,\text{in}\, \Omega_0\setminus\overline U.
$$
Give $F\in (H^1_0(\Omega_0))^*$ with $\text{supp}\,F\subset\Omega_0\setminus\overline\Omega$.
Then we get
$$\begin{array}{c}
\displaystyle
0=f(\mbox{\boldmath $G$}F\vert_U)=\tilde{f}(\mbox{\boldmath $G$}F)\\
\\
\displaystyle
=\int_{\Omega_0}\gamma_0\nabla(\mbox{\boldmath $G$}
\tilde{f})\cdot\nabla(\mbox{\boldmath $G$}F) dx\\
\\
\displaystyle
=\int_{\Omega_0}\gamma_0\nabla(\mbox{\boldmath $G$}F)\cdot\nabla(\mbox{\boldmath $G$}\tilde{f})dx\\
\\
\displaystyle
=F(\mbox{\boldmath $G$}\tilde{f}).
\end{array}
$$
Since $F$ is arbitrary, we get
$$\displaystyle
\mbox{\boldmath $G$}\tilde{f}=0\,,\text{in}\, \Omega_0\setminus\overline\Omega.
$$
From the unique continuation property(see \cite{AM} for $n=2$) and
the connectedness of $\Omega\setminus\overline U$
we get
$$\displaystyle
\mbox{\boldmath $G$}\tilde{f}=0\,\,\text{in}\,\Omega\setminus\overline U.
$$
Now let $u\vert_U\in X$ where $u\in H^1(V)$ in an open
neighbourhood $V$ of $\overline U$ and satisfies
$\nabla\cdot\gamma_0\nabla u=0$ in $V$.
By cutting off $u$ outside $\overline U$, we know that
there exists $\tilde{u}\in H^1_0(\Omega_0)$ such that
$\tilde{u}\vert_U=u\vert_U$.
Then
$$\begin{array}{c}
\displaystyle
f(u\vert_U)=f(\tilde{u}\vert_U)=\tilde{f}(\tilde{u})\\
\\
\displaystyle
=\int_{\Omega_0}\gamma_0\nabla(\mbox{\boldmath $G$}\tilde{f})\cdot\nabla\tilde{u}dx\\
\\
\displaystyle
=\int_U\gamma_0\nabla(\mbox{\boldmath $G$}\tilde{f})\cdot\nabla u dx\\
\\
\displaystyle
=\int_V\gamma_0\nabla u\cdot\nabla(\mbox{\boldmath $G$}\tilde{f})dx
=0.
\end{array}
$$
This completes the proof.

\noindent
$\Box$

\end{document}